\documentclass[reqno]{article}
\usepackage{authblk}
\usepackage{amssymb,latexsym,amsmath,amsfonts}
\usepackage{graphicx}
\usepackage{color,enumerate}
\usepackage{amsmath,amsfonts,amssymb}
\usepackage[mathscr]{eucal}
\usepackage{bbm}
\usepackage{booktabs}
\usepackage{float}
\usepackage{hyperref}
\usepackage{amsthm}
\usepackage[mathscr]{eucal}
\usepackage{paralist}
\usepackage{rotating}
\usepackage{multirow}
\usepackage[left=1.25in, right=1.25in,bottom=1.25in]{geometry}

\theoremstyle{definition}

\theoremstyle{remark}

\newcommand{\RN}[1]{\uppercase\expandafter{\romannumeral#1}}

\newcommand{\SL}[1]{\textcolor{black}{#1}}
\newcommand{\CS}[1]{\textcolor{black}{#1}}
\newcommand{\edit}[1]{\textcolor{black}{#1}}

\providecommand{\keywords}[1]{\textbf{Keywords: } #1}

\providecommand{\ams}[1]{\textbf{Mathematics Subject Classification(s): } #1}

\newcommand\extrafootertext[1]{%
	\bgroup
	\renewcommand\thefootnote{\fnsymbol{footnote}}%
	\renewcommand\thempfootnote{\fnsymbol{mpfootnote}}%
	\footnotetext[0]{#1}%
	\egroup
}

\newcommand\correspondingauthor{\thanks{Corresponding author.}}
\title{Spatial Dynamics with Heterogeneity}
\author[1]{Denis D. Patterson\correspondingauthor}
\affil[1]{High Meadows Environmental Institute, Princeton University, Princeton, NJ (denispatterson@princeton.edu)}
\author[2]{Simon A. Levin}
\affil[2]{Department of Ecology and Evolutionary Biology, Princeton University, Princeton, NJ (slevin@princeton.edu)}
\author[3]{ A. Carla Staver}
\affil[3]{Department of Ecology and Evolutionary Biology, Yale University, New Haven, CT (carla.staver@yale.edu)}
\author[4,5]{Jonathan D. Touboul}
\affil[4]{ Department of Mathematics, Brandeis University, Waltham MA  (jtouboul@brandeis.edu)}
\affil[5]{Volen Centre for Complex Systems, Brandeis University, Waltham MA}

\date{\today}
\begin{document}

\maketitle

\begin{abstract}
Spatial systems with heterogeneities are ubiquitous in nature, from precipitation, temperature and soil gradients controlling vegetation growth to morphogen gradients controlling gene expression in embryos. Such systems, generally described by nonlinear dynamical systems, often display complex  parameter dependence and exhibit bifurcations. The dynamics of heterogeneous spatially extended systems passing through bifurcations are still relatively poorly understood, yet recent theoretical studies and experimental data highlight the resulting complex behaviors and their relevance to real-world applications. We explore the consequences of spatial heterogeneities passing through bifurcations via two examples strongly motivated by applications. These model systems illustrate that studying heterogeneity-induced behaviors in spatial systems is crucial for a better understanding of ecological transitions and functional organization in brain development. 
\end{abstract} %  250 limit
\vspace*{10pt}

\vspace*{10pt}

\ams{35B32, 45K05, 92D40, 92B05}

\extrafootertext{Simon Levin and Denis Patterson thank the NSF for support via the grant DMS-1951358, and Simon Levin also appreciates support from the Army Research Office Grant W911NF-18-1-0325. Carla Staver appreciates support from NSF Grant DMS-1951394 and Jonathan Touboul appreciates support from NSF Grant DMS-1951369.}

\newpage

\section{Introduction}
In a variety of spatial natural systems, nonlinear dynamics and spatial interaction are affected by heterogeneities. Such systems have been widely studied, in a variety of settings. Systems with gradients that transition from one dominant state to another, and including a spatially transient region of bistability, generally show a \SL{more-or-less} sharp transition between the two regions; the position of this transition point depends on the intrinsic nonlinear dynamics as well as the nature of the spatial interactions. The existence of such solutions in bistable systems was recently rigorously established for a wide class of competitive reaction-diffusion equations in the limit of vanishing diffusion~\cite{perthame2015competition}; other aspects of sharp boundary formation in heterogeneous media have also been considered by numerous authors~\cite{dirr2006pinning,mori2008wave,mori2011asymptotic,xin2000front}. Systems with transient bistability are only the simplest case of what we will refer to as spatially transient bifurcation scenarios, \SL{which} we define as spatially extended systems with heterogeneities whereby the parameter varying in space is associated with the crossing of one or more  bifurcations. The behavior of more general spatially transient bifurcations remains largely open and, to the best of our knowledge, has not been addressed in a systematic manner. We argue however that understanding these dynamics is essential to better describe the behavior of complex biological systems in heterogeneous environments. We study here two particular systems paradigmatic for this problem: an ecological system describing the competition between savannas and forests within a gradient of precipitation, as well as a recently developed model of brain development featuring gradients of morphogens.

The models presented here exhibit a rich phenomenology that, in line with the topic of this special issue, have strong potential both for applications and mathematical developments. In the case of front pinning stemming from transient passage through a bistable regime, first applications of the so-called Maxwell point theory to ecology appeared recently~\cite{van2015resilience}, and were applied to the so-called Staver-Levin forest-savanna model~\cite{staver2011tree,staver_levin_2012} in the presence of a gradient in~\cite{wuyts2019tropical}, or in a \CS{similar forest-savanna} model in~\cite{goel2020dispersal}. There are also some classical results in mathematical ecology regarding the ability of species to invade, survive and coexist in heterogeneous environments~\cite{belgacem1995effects,cantrell1993permanence,cantrell1996ecological}. The recent work \cite{bastiaansen2022fragmented} on climate tipping considers a bistable nonspatial model subject to diffusion and a non-monotonic heterogeneous forcing term, and is arguably closest to the mechanisms studied here. Similarly, researchers in mathematical neuroscience have studied how the shape of the input into neural field models can induce interesting dynamics in the output function~\cite{faye2014pulsatile,kilpatrick2008traveling,kilpatrick2013optimizing,kilpatrick2013wandering}. However, in these cases, the complexity of the dynamics are induced by the shape of the gradients. In contrast, we will focus on simpler monotonically varying heterogeneities and, as we will show, complex dynamics can emerge simply from the interaction of the underlying non-spatial dynamics and the slope of the gradient (or the length scale of the problem). Moreover, we aim to highlight the importance of the speed (or size of the region of space) for which the heterogeneity crosses different dynamical regimes of the underlying non spatial model. As we illustrate, this crucially influences the emergence or non emergence of new spatial patterns or spatio-temporal dynamics.

Turing pattern formation in heterogeneous systems (but not including transient passages through Turing instabilities) has also been widely explored, and it was shown that dispersal can sharpen boundaries~\cite{garcia2000dispersal}. Numerous authors have developed criteria for pattern formation in specific reaction-diffusion models where some parameters are allowed to vary across the spatial domain in a stepwise manner~\cite{benson1993diffusion,page2003pattern}. However, we believe that smooth variation of the gradient to be a key ingredient in many applications (cf. section \ref{sec.brain_patterning}); it changes the mathematical approaches to the problem of characterizing patterns significantly and also produces qualitatively different solutions, e.g. blending of spot, stripes and labyrinths, as well as multi-frequency patterns~\cite{feng2021coup}. Researchers have also pursued perturbative results about the homogeneous case. For example, Benson et al.~\cite{benson1998unravelling} carried out a remarkably detailed bifurcation analysis of a reaction-diffusion system with diffusion coefficient $\mathcal{D}(x) = D + \eta x^2$ in the limit as $\eta \downarrow 0$, illustrating that a rich array of dynamics can be spawned by small heterogeneity. The dynamics of spiked solutions to pattern forming reaction diffusion systems subject to spatially localized heterogeneity have also attracted considerable attention (see \cite{wong2021spot} and the references therein), while other work on this topic has allowed more general heterogeneous structure by working in asymptotic regimes of the diffusivity strengths~\cite{ward2002dynamics}. However, biological problems (such as the one we present in section \ref{sec.brain_patterning}) often exhibit spatial heterogeneity on the same scale of the spatial domain itself, rendering many asymptotic or perturbative results of relatively limited utility in applications. 

There has been recent and growing interest in revisiting and generalizing the classical Turing pattern paradigm to account for general spatial heterogeneity~\cite{krause2021modern,van2021pattern}. Although the literature has focused mainly on heterogeneity in the diffusion coefficient(s), several investigators have recently employed series-approximation approaches to derive conditions for the onset of patterns with heterogeneity in the reaction dynamics~\cite{kozak2019pattern,krause2020one}. While these results are undoubtedly promising, the vast bulk of the existing literature deals with reaction-diffusion systems. Many applications require models incorporating more complex spatial operators, such as chemotaxis, nonlocal spatial operators or advection (which can alter pattern forming dynamics even when present at very low levels~\cite{klika2018domain}), thus demanding the development of new mathematical machinery. The inclusion of heterogeneity in pattern forming systems has also recently been reported by several authors to induce spatio-temporal oscillations~\cite{dornelas2021landscape,krause2018heterogeneity}, posing yet more interesting mathematical questions in this area. 

In section \ref{sec.rainfall_gradient}, we begin our exploration of heterogeneity-induced dynamics by first investigating the forest-grass Staver-Levin model with nonlocal interactions and a rainfall gradient that increasingly favors forest tree expansion. Turning next to a more detailed version of the model with four functional types, we incorporate the impact of the rainfall on the various types of vegetation by allowing the forest and savanna tree birth rates to vary spatially, since they are likely to exhibit the most spatial variation and impact outcomes most strongly. Our approach is to adopt qualitatively appropriate gradients, as opposed to trying to quantitatively estimate the impact of increasing rainfall on the productivity of each of the functional types in the model. We are guided by a detailed bifurcation analysis of the corresponding nonspatial model and this allows us to predict potential emergent behaviors caused by the heterogeneity and illustrate interesting spatio-temporal dynamics. Through a numerical bifurcation analysis of the system of integro-differential equations, we find multistability between several nonhomogeneous solutions for a significant range of the dispersal parameter; we find \SL{front-pinned} solutions (which connect stable solution branches from the nonspatial model), and nonhomogeneous grass and forest dominated solutions. In the \CS{four-functional-type} spatial Staver-Levin model there is a much greater variety of solutions in the presence of heterogeneity. We highlight several plausible transition mechanisms for the \SL{empirically} observed savanna to forest transition with increasing rainfall. Once more, we uncover \SL{front-pinned} solutions, including \SL{a front-pinned solution predicting} a grassland band mediating the savanna-forest transition, as well as multistability and changes in stability as the dispersal parameters vary. Additionally, we observe periodic waves of invasion for gradients \SL{that} intersect a large region of parameter space that produces stable oscillations in the nonspatial Staver-Levin model. Intriguingly, as we vary the width or speed at which the gradient cuts through oscillatory region, we see period doubling of the waves and eventually more complex solutions \SL{that} appear to display the hallmarks of spatio-temporal chaos. This illustrates how the speed at which the gradient crosses through different dynamical regimes or bifurcations crucially determines the resulting dynamics.

In section \ref{sec.brain_patterning}, we discuss a PDE model of arealization in the mammalian forebrain featuring transcription factor gradients that mediate the competition between different neural fates in early embryonic development~\cite{feng2021coup}. This model explained unexpected patterning and dislocation of the typically sharp boundary between abutting cortical regions upon artificial manipulation of the transcription factor gradients in mice. Mathematically, the dramatic breakdown of the \SL{boundary-forming} mechanism in the mice corresponds to a transient passage of the gradient through a region of parameter space with \SL{pattern-forming} instabilities. For a 1D spatial domain, we demonstrate that the structure of the resulting solutions of the mathematical model depend crucially on how the gradient transits the pattern forming region and emphasize the qualitative differences when comparing these solutions to those for a typical pattern forming system on a homogeneous domain.

\section{The savanna-forest transition with a rainfall gradient}\label{sec.rainfall_gradient}
The first example we shall study is a transient passage through a family of periodic orbits, in the case of the transition between savannas and forest within a gradient of precipitation. 
\subsection{Spatially extended savanna-forest dynamics }
The vegetation model that we shall use for this, hereafter referred to \SL{as the} Staver-Levin (SL) model~\cite{staver_levin_2012}, describes the interaction between savanna trees with an age structure ($S$ for savanna saplings, $T$ for adult savanna trees), forest trees ($F$), and grass patches ($G$). Grass  patches are locations that carry fires that limit the expansion of forest trees and delay the maturation of savanna saplings, but are also the locations where new trees of both types can grow. In~\cite{patterson2019probabilistic}, we introduced a spatially explicit stochastic model accounting for birth, death and interaction rules between these ecological species proved that the fractions of the different components at location $x$ and time $t$\footnote{Also representing the probability for a site at location $x$ to be of a given type at time $t$.} satisfy the following system of nonlinear integro-differential equations in an appropriate mean-field limit:
\begin{subequations}\label{SL_integral}
	\begin{alignat}{2}
		\partial_t G(x,t) &= \mu S + \nu T + \phi\left(\int_\Omega w(x-y) G(y,t)\,dy \right)F  - \alpha G \int_\Omega J_F(x-y)F(y,t)\,dy \nonumber\\
		&\quad - \beta \,G \int_\Omega J_T(x-y)T(y,t)\,dy, \\
		\partial_t S(x,t) &= -\mu S - \omega\left(\int_\Omega w(x-y)G(y,t)\,dy \right) S - \alpha S \int_\Omega J_F(x-y)F(y,t)\,dy \nonumber\\ &\quad + \beta \,G  \int_\Omega J_T(x-y)T(y,t)\,dy,\\
		\partial_t T(x,t) &= -\nu T + \omega\left( \int_\Omega w(x-y)G(y,t)\,dy \right) S - \alpha T  \int_\Omega J_F(x-y)F(y,t)\,dy, \\
		\partial_t F(x,t) &= \alpha (G + S + T) \int_\Omega J_F(x-y)F(y,t)\,dy - \phi\left(\int_\Omega w(x-y) G(y,t)\,dy\right)F,
	\end{alignat}
\end{subequations}
for each $(x,t) \in \Omega \times \mathbb{R}^+$ for some $\Omega \subset \mathbb{R}^2$.

Since the system of equations given by \eqref{SL_integral} describes the evolution of probability densities, we have
\begin{equation}\label{eq.normalization}
	G(x,t) + S(x,t) + T(x,t) + F(x,t) = 1 \mbox{ for each } (x,t)\in \Omega\times \mathbb{R}^+.
\end{equation}
In this model, the constants $\mu$ and $\nu$ are the mortality rates of saplings and savanna trees respectively. The functions $\phi$ and $\omega$ represent the burning rates of forest trees and saplings due to fire; they are theoretically predicted, and empirically observed, to have sharp threshold or sigmoidal profiles as \SL{functions} of the available flammable cover (grass in this framework)~\cite{schertzer2015implications}. For numerical investigations, we employ the following smooth approximations to a sigmoid for $\phi$ and $\omega$:
\begin{equation}\label{eq.sigmoids}
\omega(G) = \omega_0 + \frac{\omega_1-\omega_0}{1 + e^{-(G-\theta_1)/s_1}}, \quad \phi(G) = \phi_0 + \frac{\phi_1-\phi_0}{1 + e^{-(G-\theta_2)/s_2}} \quad \mbox{for } G\in[0,1].
\end{equation}
with parameter values as given in Table \ref{table.SL_parameters} below.

The kernel function $w$ measures the ability of fire to spread spatially from a point \SL{that} is already burning. The constants $\alpha$ and $\beta$ account for the strength of \SL{forest-tree} and \SL{savanna-tree} invasion via seed dispersal, with the spatial distribution of these seeds captured by the kernels $J_F$ and $J_T$. The inclusion of nonlocal or long-range interactions is considered by many most appropriate for spatial vegetation models as dispersal of seeds is often long range or even \SL{heavy-tailed}~\cite{nathan2012dispersal,thompson2008plant}. The spatial interaction (fire spread and seed dispersal) are assumed isotropic so all kernels are of convolution type and the model \eqref{SL_integral} is thus (for now) posed on a homogeneous spatial domain. For simplicity, and in all numerical results, we use zero mean Gaussian kernels with different standard deviations (to reflect the relative length scales of the different spatial process). In particular, we have
\[
\mathcal{G}(x,\sigma) := \frac{1}{\sqrt{2\pi \sigma^2}}e^{-x^2/2\sigma^2}, \quad\sigma >0,\quad   x \in \Omega,
\]
with $w(x) = \mathcal{G}(x,\sigma_W)$, $J_T(x) = \mathcal{G}(x,\sigma_T)$, and $J_F(x) = \mathcal{G}(x,\sigma_F)$. In the main text, we will always consider the system with \emph{reflecting boundary conditions}, which are appropriate for both the nonlocal operators and the heterogeneous nature of the models (see Appendix \ref{sec.boundary_conditions}). The results are qualitatively similar for open boundary conditions and the Supplementary Materials contain simulations with open boundaries for comparison purposes.

In this work, we allow $\alpha$ and $\beta$ to be functions of spatial position, $x$, to mimic heterogeneous effects due to environmental gradients observed in reality, particularly the pronounced rainfall gradients in sub-Saharan Africa and the Amazon~\cite{bucini2007continental,wuyts2017amazonian}. We will neglect the influence of environmental variation on the fire and natural mortality processes in the model since we expect these effects to be less impactful on the dynamics. \CS{However, we note that fire frequency and intensity tend to increase and then decrease with increasing rainfall~\cite{he2018baptism} and that tree mortality tends to decrease and then increase with increasing rainfall and varies substantially with soil properties~\cite{quesada2012basin}.}

Table \ref{table.SL_parameters} below summarizes the parameters of the SL model, along with their ecological interpretations and default numerical values. Further details on the parameter gradients and numerical schemes can be found in Appendix \ref{sec.numerical_parameters}.
\begin{table}[H]
	\centering
	\caption{Summary of parameters for the SL model}
	\begin{tabular}{@{}lcc@{}}
		\toprule
		Ecological interpretation  & Expression & Default Value \\
		\midrule
		Forest tree birth rate & $\alpha(x)$ & [piecewise linear] \\   
		Savanna saplings birth rate & $\beta(x)$ & [piecewise linear]  \\   
		Savanna sapling-to-adult recruitment rate & $\omega$ & $\omega(G) = \omega_0 + \frac{\omega_1-\omega_0}{1 + e^{-(G-\theta_1)/s_1}}$ \\   
		 &  & $\omega_0=0.9$, $\omega_1 = 0.4$,  \\ 
		 & & $\theta_1 = 0.4$, $s_1 = 0.01$ \\
		Forest tree mortality rate  & $\phi$ & $\phi(G) = \phi_0 + \frac{\phi_1-\phi_0}{1 + e^{-(G-\theta_2)/s_2}} $  \\
		 &  & $\phi_0=0.1$, $\phi_1 = 0.9$,   \\
		 & & $\theta_2 = 0.4$, $s_2 = 0.05$ \\
		Savanna sapling mortality rate & $\mu$ & 0.1  \\
		Adult savanna tree mortality rate & $\nu$ & 0.05 \\
		Forest tree seed dispersal parameter & $\sigma_F$ & -  \\
		Savanna tree seed dispersal parameter & $\sigma_T$ & - \\
		Fire spread (dispersal) parameter & $\sigma_W$ & - \\
		\bottomrule
	\end{tabular}
	\label{table.SL_parameters}
\end{table}

In addition to the \CS{four-functional-type} SL model stated above, numerous other related models with the same (or very similar) underlying interaction rules have been studied in the  literature~\cite{durrett2015coexistence,durrett2018heterogeneous,hoyer2021impulsive,li2019spatial,schertzer2015implications,touboul2018complex,wuyts2019tropical}. One important special case is that in which only the grass and forest types are present, essentially reducing the model to a direct competition between fire and seed dispersal effects. In this case, allowing $\alpha$ to vary in space, the system \eqref{SL_integral} reduces to the single integro-differential equation
\begin{multline}\label{eq.GrassForest_space}
	\partial_t G(x,t) = \phi\left( \int_{\Omega}w(x-y)G(y,t)\,dy \right)\left(1 - G(x,t)\right) - \alpha(x) G(x,t)\left(1 - \int_{\Omega}J_F(x-y)G(y,t)\,dy\right), 
\end{multline}
where the normalization condition $G(x,t)+ F(x,t) = 1$ and the assumption that $J_F$ and $w$ are probability density functions on $\Omega$ (and hence have unit integrals) enable this reduction.

\subsection{Heterogeneities}
In the tropics, empirical data \SL{show} that while savanna landscapes tend to dominate at low rainfall levels, there is evidence of bistability between \SL{savanna-dominated} and \SL{forest-dominated}  biomes at intermediate rainfall, and forest naturally tends to dominate at sufficiently high rainfall levels \cite{staal2020hysteresis,staver2011tree}. Consequently, any model purporting to explain the savanna-forest transition should explicitly account for the spatial heterogeneity arising from the significant differences in rainfall levels across the domain. To this end, and for the sake of tractability, we consider a one-dimensional bounded spatial domain $\Omega = [0,1]$ throughout and allow both the \SL{forest-tree} birth rate $\alpha$ and the \SL{savanna-tree} birth rate $\beta$ to vary with position in the domain; in particular, we choose both to be increasing linear functions of $x\in\Omega$ so that the domain becomes progressively wetter (and more forest favored) as $x$ increases~(see \cite{staver_levin_2012} and the references therein). Intuitively, we expect the slope of the forest birth rate to be larger than that of the savanna trees since eventually the forest trees must \SL{out-compete} the savanna trees at high rainfall levels.

We primarily study the dynamics of the heterogeneous model in the \SL{small-dispersal} limit since the seed and fire dispersal scale is several orders of magnitude smaller than the scale on which the rainfall gradient is significant. We focus on the case of centered Gaussian dispersal kernels as a representative class of kernels which illustrate our main qualitative conclusions. However, our conclusions should hold for a broad class of nonnegative even kernels, although spreading speeds will differ for heavy-tailed kernels and it would undoubtedly be interesting to investigate the impact of heavy-tailed dispersal on the nonequilibrium solutions shown in section \ref{sec.waves}.
	
\subsection{Front pinning in the grass-forest model}\label{sec.rainfall_GF}
\begin{figure}[h]
	\centering
	\includegraphics[width=0.9\textwidth]{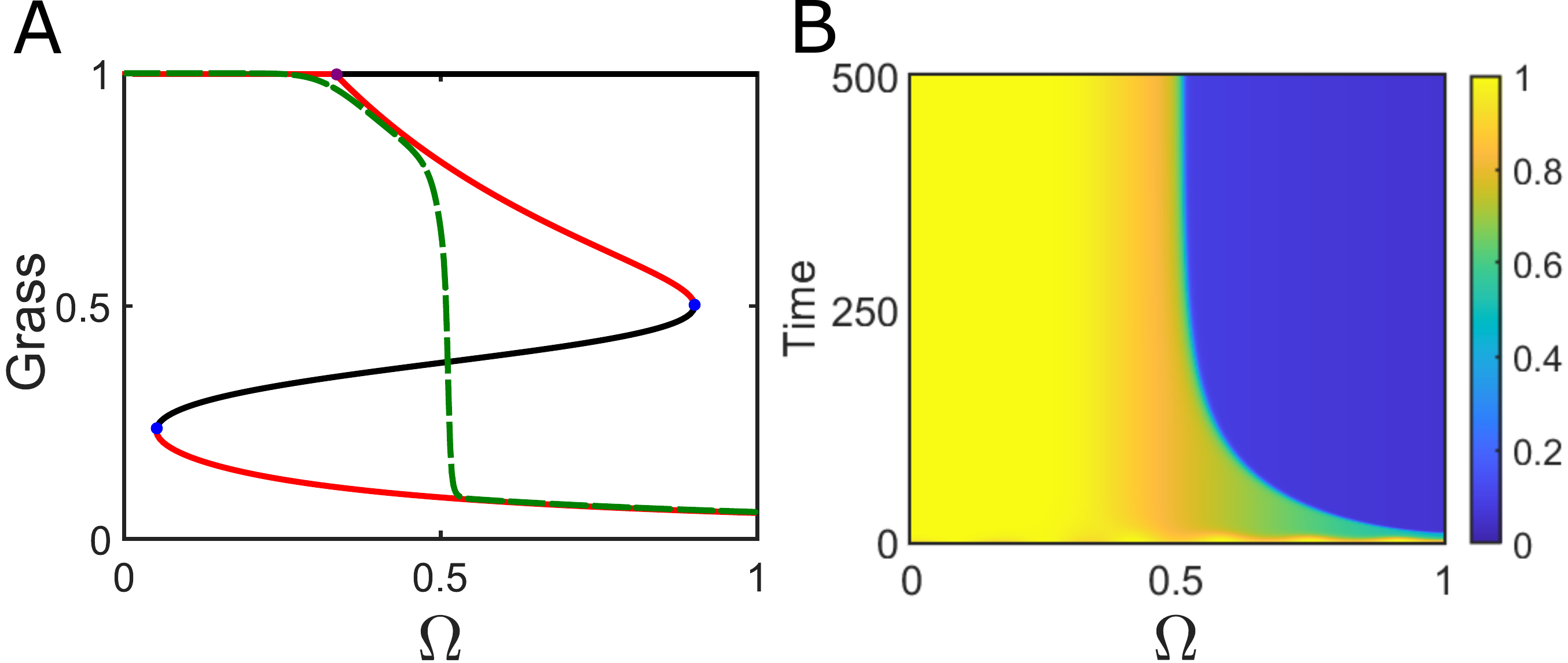}
	\caption{\textbf{A:} ODE bifurcation diagram with space as the bifurcation parameter (stable equilibrium solution curves in red and unstable equilibria in black) overlaid with a stable steady-state solution of the IDE for dispersal $\sigma = 0.01$ on the heterogeneous domain (dashed green line). \textbf{B:} Space time plot of convergence to the steady-state solution plotted in A; the speed of the wave of invasion forest into grass slowly tends to zero before equilibrium is reached.}\label{fig.front_pinning_1}
\end{figure}

Consider the grass-forest model given by equation \eqref{eq.GrassForest_space}  on the 1D spatial domain $\Omega = [0,1]$. If we assume completely localized interactions, our model reduces to a family of ordinary differential equations indexed by the spatial variable $x$:
\begin{equation}\label{eq.GrassForest_ODE}
	\dot{G} = \phi(G)(1-G) - \alpha(x)G(1-G),\quad x\in [0,1],
\end{equation}
Equation \eqref{eq.GrassForest_ODE} thus describes the evolution of the grass cover proportion in an isolated patch of landscape (i.e. one without spatial interactions). We study the corresponding \SL{spatial} model, given by \eqref{eq.GrassForest_space}, subject to a rainfall gradient which spans from a \SL{grass-dominated} regime at the left \SL{end-point} of the domain $x=0$ to a \SL{forest-dominated} regime at the \SL{right-most} point of the domain $x=1$; we reflected this in the model by choosing $\alpha(x)$ linear with positive slope. Such a gradient in the forest tree birth rate reflects the increasing productivity of forest trees at progressively higher rainfall. Figure \ref{fig.front_pinning_1} B shows the \SL{time-evolution} of grass cover in the spatial grass-forest model \eqref{eq.GrassForest_space} with dispersal parameter $\sigma=0.01$ with rainfall gradient chosen as outlined above. Figure \ref{fig.front_pinning_1} A shows the final steady-state solution profile from Figure \ref{fig.front_pinning_1} B (dashed green line) overlaid with the bifurcation diagram from the corresponding nonspatial model \eqref{eq.GrassForest_ODE} (red for stable equilibria, black for unstable). 

Solutions to equation \eqref{eq.GrassForest_ODE} and their stability are computed for each value of $x$ in XPP using AUTO~\cite{ermentrout2003simulating}. Solutions to the corresponding spatial model \eqref{eq.GrassForest_space} are computed using a finite-difference spatial discretization of the integrodifference equation so that it can be solved via time stepping schemes and with a fixed point approach in MATLAB for bifurcation analysis.  We recall that all simulations shown in the main text are with Gaussian kernels and for reflecting boundary conditions (see Appendix \ref{sec.boundary_conditions}). However, the results are very similar for other centered univariate kernels and for other types of boundary conditions since in the small to moderate dispersal cases we  study, boundary effects are minimal. We also systematically explored the space of initial conditions to uncover multistability, which was observed in a number of cases. Grid convergence and time-step convergence studies were carried out to verify that the stability of computed solutions was robust.

By overlaying the nonspatial bifurcation diagram onto the spatial solution in Figure \ref{fig.front_pinning_1} A, we see that the spatial solution essentially interpolates between the two stable solution branches from the nonspatial model with a sharp transition between the grass and forest dominated parts of the domain. As per Figure \ref{fig.front_pinning_1} A, the gradient begins before the first saddle node bifurcation of the nonspatial model and hence the lower solution branch is not a candidate solution across the whole domain in the spatial model, as we might have expected. Likewise, the gradient extends beyond the upper nonspatial \SL{saddle-node} bifurcation and so the upper solution branch is similarly not a candidate spatial solution either (at least for small dispersal). This is a completely new type of solution not present without spatial interactions and heterogeneity. Numerous authors have proposed this type of mechanism as an explanation for the empirically observed sharp savanna-forest boundaries in sub-Saraharan Africa, the Amazon, and other tropical regions with similarly abrupt savanna-forest transitions~\cite{champneys2021bistability,goel2020dispersal,goel2020dispersal_madagascar,wuyts2017amazonian}.

Convergence to the steady-state in Figure \ref{fig.front_pinning_1} B proceeds via a wave of invasion that fails to propagate beyond a point in the spatial domain and where a sharp front forms in the solution. This phenomenon has been termed ``front pinning'' or ``range pinning'' in the literature and has been demonstrated in many bistable \SL{PDE-based} ecological models~\cite{champneys2021bistability,van2015resilience,wuyts2019tropical}, and other applied contexts~\cite{kulka1995influence}. The point in the spatial domain at which the front forms in the solution is referred to as the \emph{Maxwell point}, and is a \edit{generically stable configuration in heterogeneous systems. Front pinned solutions can also be stable in certain spatially homogeneous systems~\cite{mori2011asymptotic,mori2008wave}.} Heuristically, \edit{pinning} is \edit{harder to observe in a homogeneous system} because waves of invasion will select some constant wave speed $c$ with which to propagate \edit{on} a homogeneous \edit{domain.} If the system parameters are set to exactly generate $c=0$, then \SL{front-pinning} can occur, but this parameter set will have measure zero in parameter space and hence a pinned front solution is not generic (and thus physically unrealistic). In our model, we have nonconstant wave speeds $c(x)$ due to the heterogeneity in the system and bistability persists for a large parameter region in the absence of spatial interactions (see Figure \ref{fig.front_pinning_1}), so stable \SL{front-pinned} solutions are expected for appropriate gradient choices. Various authors have considered calculation of the Maxwell point for reaction diffusion versions of the Staver-Levin model and other ecological models~\cite{wuyts2019tropical,goel2020dispersal}. However, this approach appears technically limited to scalar systems with diffusion as the dominant spatial interaction, meaning the method is approximate for models with nonlocal interactions. Moreover, the Maxwell point calculations are only approximate when the values of the nonspatial steady states vary as a function of the gradient (as they do in our model, see Figure \ref{fig.front_pinning_1} A). The existence of \SL{front-pinned} solutions has been established rigorously in reaction-diffusion models via asymptotic analysis in the small diffusion limit~\cite{perthame2015competition}, but, to the best of our knowledge, this question remains open for systems involving nonlocal operators.

\begin{figure}[h]
	\centering
	\includegraphics[width=\textwidth]{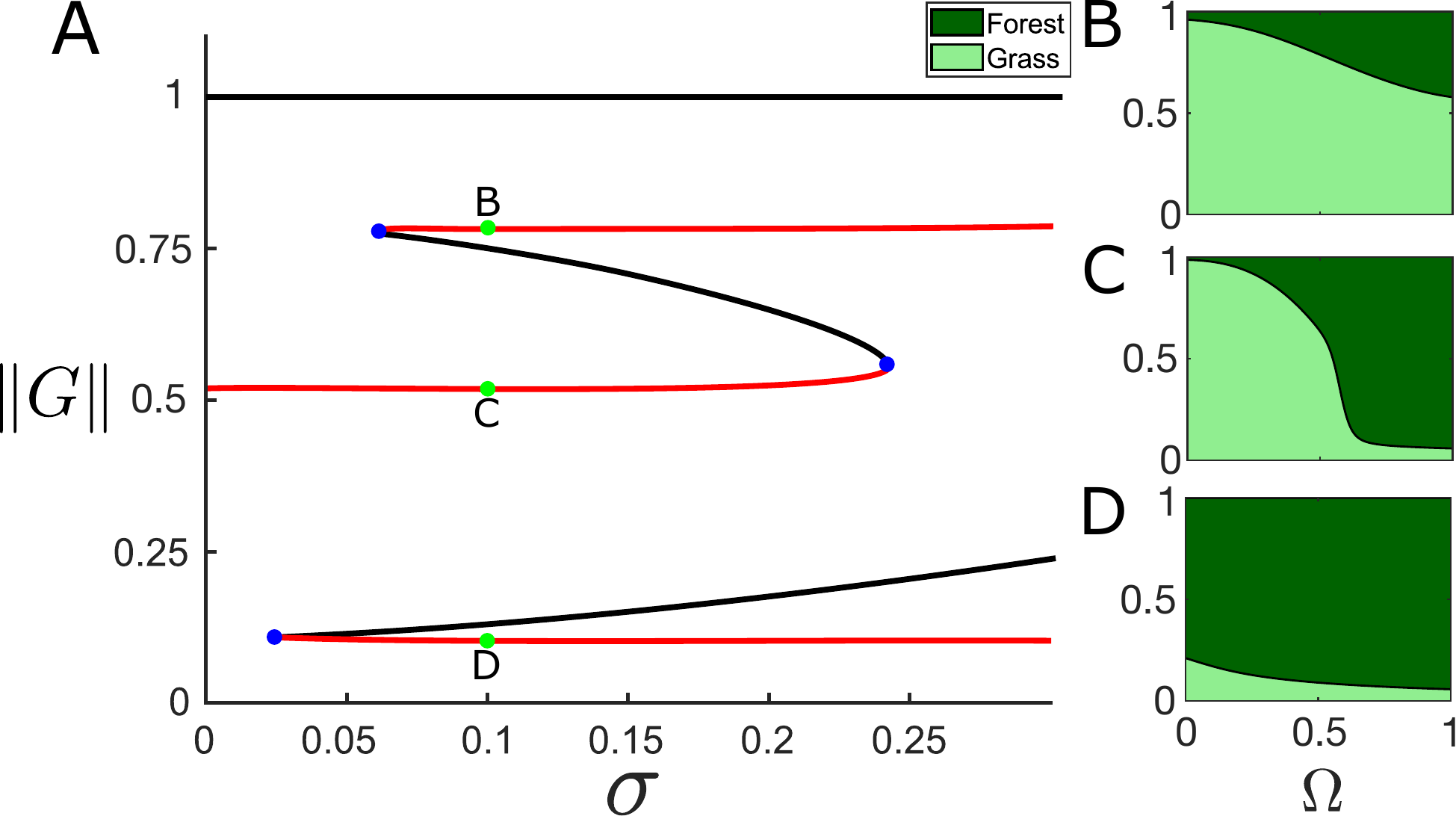}
	\caption{\edit{\textbf{A:} Bifurcation diagram as the dispersal parameter $\sigma= \sigma_F = \sigma_W$ is varied with the $L^1$ norm of the grass component of the solution on the $y$-axis (reflecting boundary conditions and Gaussian kernels for both seeds and fire). Red lines: Stable solutions, Black lines: Unstable solutions,  \textbf{B/C/D:} Solutions from the grass dominated branch, the front-pinned branch and the forest branch, respectively.}}\label{fig.bifurcations_dispersal}
\end{figure}
% notes on numerical parameters: N=400 grid points, fsolve tolerance 1e-10, L^1 tolerance 1e-8 (most solutions have L^1 error of around 1e-16)

The existence and stability of these solutions are functions of the spatial interactions.  \edit{In Figure \ref{fig.bifurcations_dispersal} we investigate the impact of the dispersal parameter $\sigma$ on front pinning and multistability in the system\footnote{\edit{Note that varying the dispersal parameter in this way is equivalent to varying the size of the spatial domain (with a smaller dispersal parameter corresponding to a larger domain).}}. The standard deviations of the Gaussian fire and seed kernels, $\sigma_W$ and $\sigma_F$, are chosen equal to a common value $\sigma$, which serves as the bifurcation parameter.} We observe that there is a large dispersal range (or range of spatial scales) for which the model supports three stable solutions, a \SL{grass-dominated} solution (Figure \ref{fig.bifurcations_dispersal} B1), a \SL{forest-dominated} solution (Figure \ref{fig.bifurcations_dispersal} B3) and a front-pinned solution with a sharp transition between forest and grass (Figure \ref{fig.bifurcations_dispersal} B2). This latter solution disappears for high dispersal values through a saddle-node bifurcation. Heuristically, when dispersal is large compared to the spatial heterogeneity, it effectively rapidly homogenizes the system and the system can \SL{no longer} support such a heterogeneous solution and instead converges to more well-mixed solutions. \edit{We further observe that the forest-dominated solution branch exists and remains stable for all but very low dispersal values, eventually disappearing in a saddle-node bifurcation around $\sigma\approx 0.026$. The grass-dominated solutions disappear at low dispersal rates, in favor of mixed equilibria associated with unavoidable forests in most humid regions and grass at driest regions. The all-grass solution (i.e. $G(x) = 1$ for all $x \in [0,1]$) is a solution for every value of $\sigma$ but is always unstable.}

\subsection{\SL{Front-pinning} in the \CS{four-functional-type} forest-savanna model}
\begin{figure}[!ht]
	\centering
	\includegraphics[width=\textwidth]{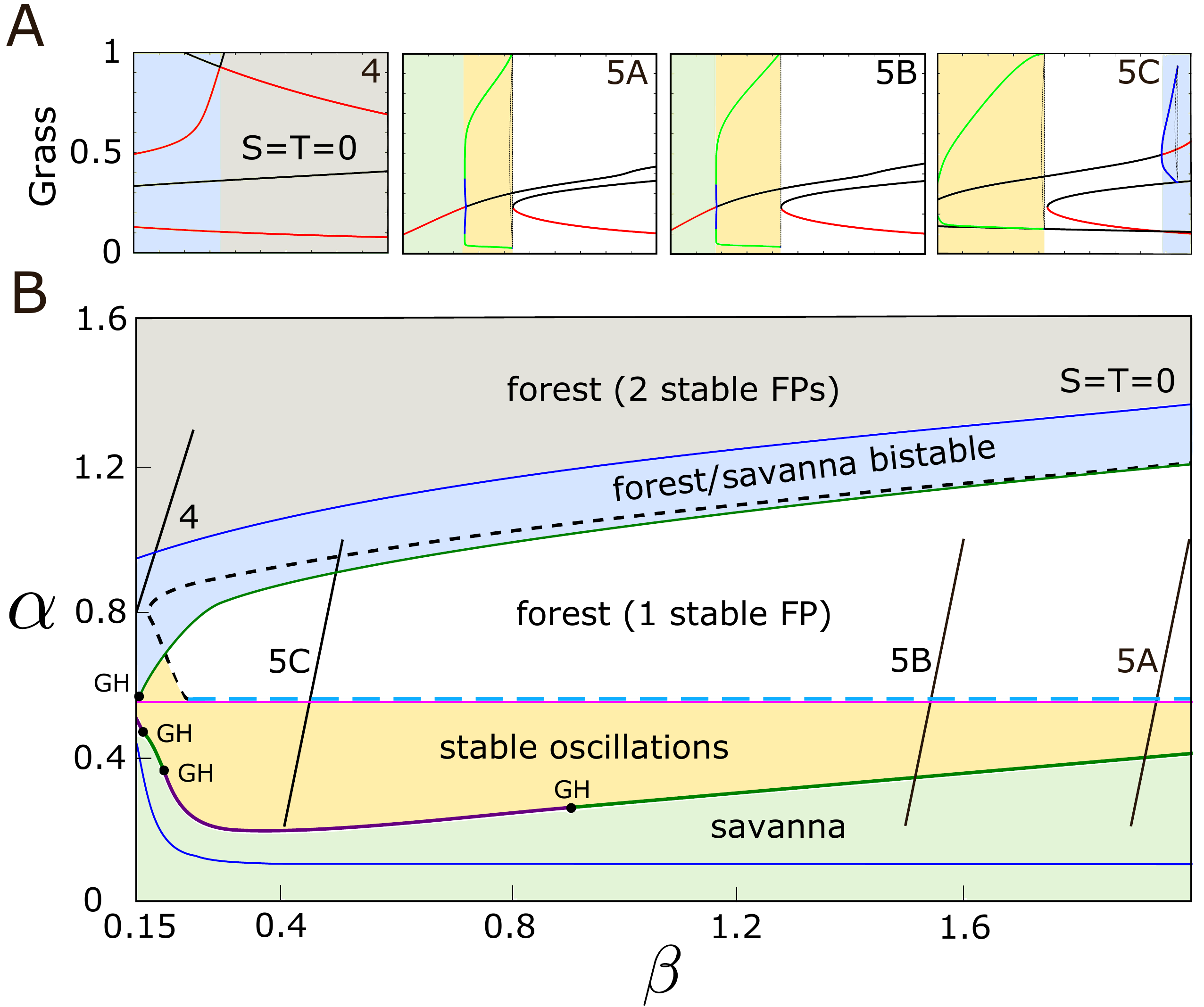}
	\caption{\edit{\textbf{A: }One-parameter bifurcation diagrams corresponding to each of the rainfall gradients studied. Stable/unstable equilibria denoted by red/black, stable/unstable limit cycles in green/blue and heteroclinics denoted by dashed black lines. \textbf{B:} Two-parameter bifurcation diagram in $\alpha$ and $\beta$ for the nonspatial Staver-Levin model (i.e. completely localized interactions). Transcritical bifurcation curves in blue, saddle node curves in magenta, supercritical Hopf curves in purple and subcritical Hopf curves in dark green, with switching points from one type to the other at Bautin (Generalized Hopf, label GH) points.}}\label{fig.codim2}
\end{figure}
We introduce the effects of a rainfall gradient phenomenologically by allowing both $\alpha$ and $\beta$ to be linearly increasing in $x$, meaning that both kinds of trees become more productive as we move along the gradient, or \CS{moving} rightwards in the 1D spatial domain $\Omega = [0,1]$. To reflect the empirical \CS{observation} that forest trees dominate at high rainfall, the slope of $\alpha$ will be larger than that of $\beta$ in all scenarios. \CS{Estimating the real quantitative impact of sub-Saharan African rainfall on each of these functional types is a formidable and outstanding challenge}, but we can glean considerable insight into the range of dynamics and savanna-forest transitions with our more qualitative approach. Previous studies have noted the array of behaviors that emerge upon introducing a realistic rainfall gradient into this model~\cite{wuyts2019tropical}, but it can be difficult to motivate and understand the resulting spatial dynamics without reference to the already complex dynamics of the nonspatial version of \eqref{SL_integral}. Figure \ref{fig.codim2} B shows the two-parameter bifurcation diagram for the (nonspatial) \CS{four-functional-type} Staver-Levin model as a function of the forest tree birth rate $\alpha$ and the savanna tree birth rate $\beta$. We immediately see that a rainfall gradient tracing a linear path in $\alpha$-$\beta$ space may connect regions with very different dynamics in the absence of spatial interactions, including stable limit cycles involving all four \CS{functional types} (see \cite{touboul2018complex} for more details). The top row of subfigures in Figure \ref{fig.codim2} A shows various one parameter bifurcation diagrams (only the grass component of the solution is shown) along gradients in $\alpha$-$\beta$ space for which we discuss the corresponding spatial dynamics below.

\begin{figure}[h]
	\centering
	\includegraphics[width=\textwidth]{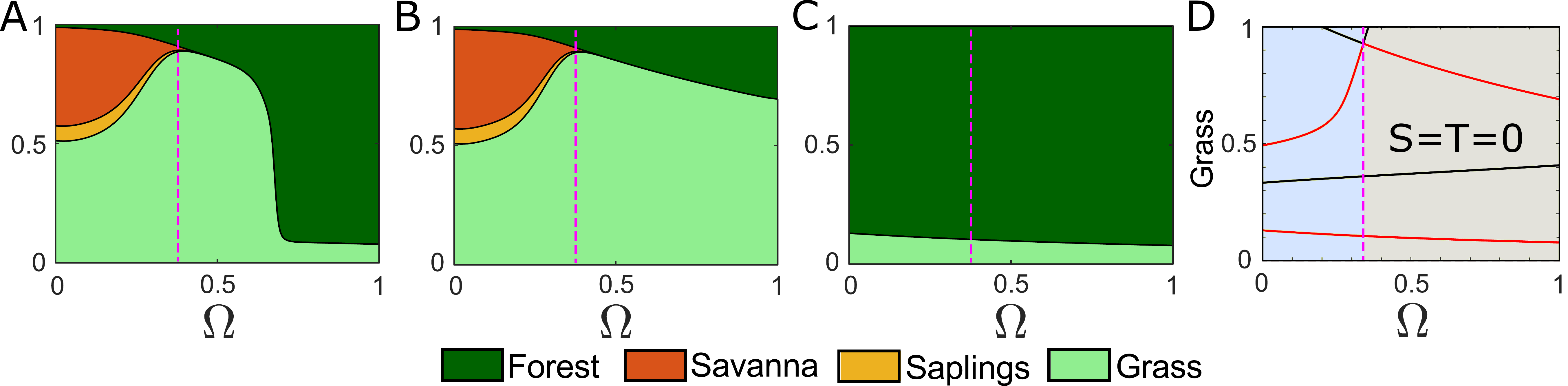}
	\caption{\textbf{A:} For $x\in[0,\,0.4]$, the system is in a savanna state, for $x\in[0.4,\,0.7]$ there is a grass dominated state with a small number of forest trees and then after the Maxwell point around $x=0.7$, the system transitions to a forest dominated state. \textbf{B:} Alternative stable state with savanna on the left side of the domain with grass dominant on the right of the domain. \textbf{C:} Alternative stable state with forest dominant across the entire domain. Dispersal parameters are set to $\sigma=0.025$ for all panels and boundary conditions are reflecting. \textbf{D:} Nonspatial bifurcation diagram for \eqref{SL_integral} with the spatial variable $x$ as the bifurcation parameter.}\label{fig.grass_front}
\end{figure}
% parameters for A and B:
% alpha = 1.05;
% alpha_s = 0.25 
% beta = 0.2
% beta_s = 0.05
% sigma = 0.02

Figure \ref{fig.grass_front} shows three multi-stable solutions for a fixed parameter set and a linear rainfall gradient (shown in $\alpha$-$\beta$ space in Figure \ref{fig.codim2}) in the savanna-forest model given by \eqref{SL_integral}; the dispersal parameters are identically $\sigma=0.025$ in all panels. Panels A and B of Figure \ref{fig.grass_front} show a pair of stable solutions in which all four \CS{functional types} are present. The solution in \ref{fig.grass_front} A follows the upper stable branch in the nonspatial bifurcation diagram (Figure \ref{fig.grass_front} D) until a Maxwell point around $x\approx 0.7$ where there is a rapid transition to forest dominance (and the lower stable branch in Figure \ref{fig.grass_front}). This solution has both a stable savanna and a stable forest domain separated by a grass band sufficiently high to suppress forest trees; this is a novel prediction in terms of the savanna-forest transition in the tropics and does not seem to appear in the \CS{existing} literature, \CS{but anecdotally is consistent with the low tree biomass and density observed in, e.g., the Bateke Plateau bordering the Congo rainforest~\cite{nieto2018mesic}}. The grass band is caused by a transcritical bifurcation in the nonspatial model (marked by the vertical pink dashed line), up to which, grass steadily increases at the expense of savanna trees; this high level of grass exceeds the ignition threshold for fire (at least locally) and hence keeps forest at a low level. At the transcritcal, the savanna loses stability and the system enters the forest-grass subsystem, in the grass-dominated state. Figure  \ref{fig.grass_front} panel C shows a forest dominated solution that is stable for all levels of the dispersal parameter. As in the forest-grass  front-pinning example, the \SL{front-pinned} solutions eventually lose stability at sufficiently large levels of dispersal (equivalent to a smaller spatial domain); this occurs around $\sigma\approx 0.05$ for the gradient chosen in this particular example. \edit{In the Supplementary Materials, we show that this example is robust to stochasticity and nonlinearity in the rainfall gradient (see SM1.).}

\subsection{Transient passage through a family of periodic orbits}\label{sec.waves}

\begin{figure}[h]
	\includegraphics[width=\textwidth]{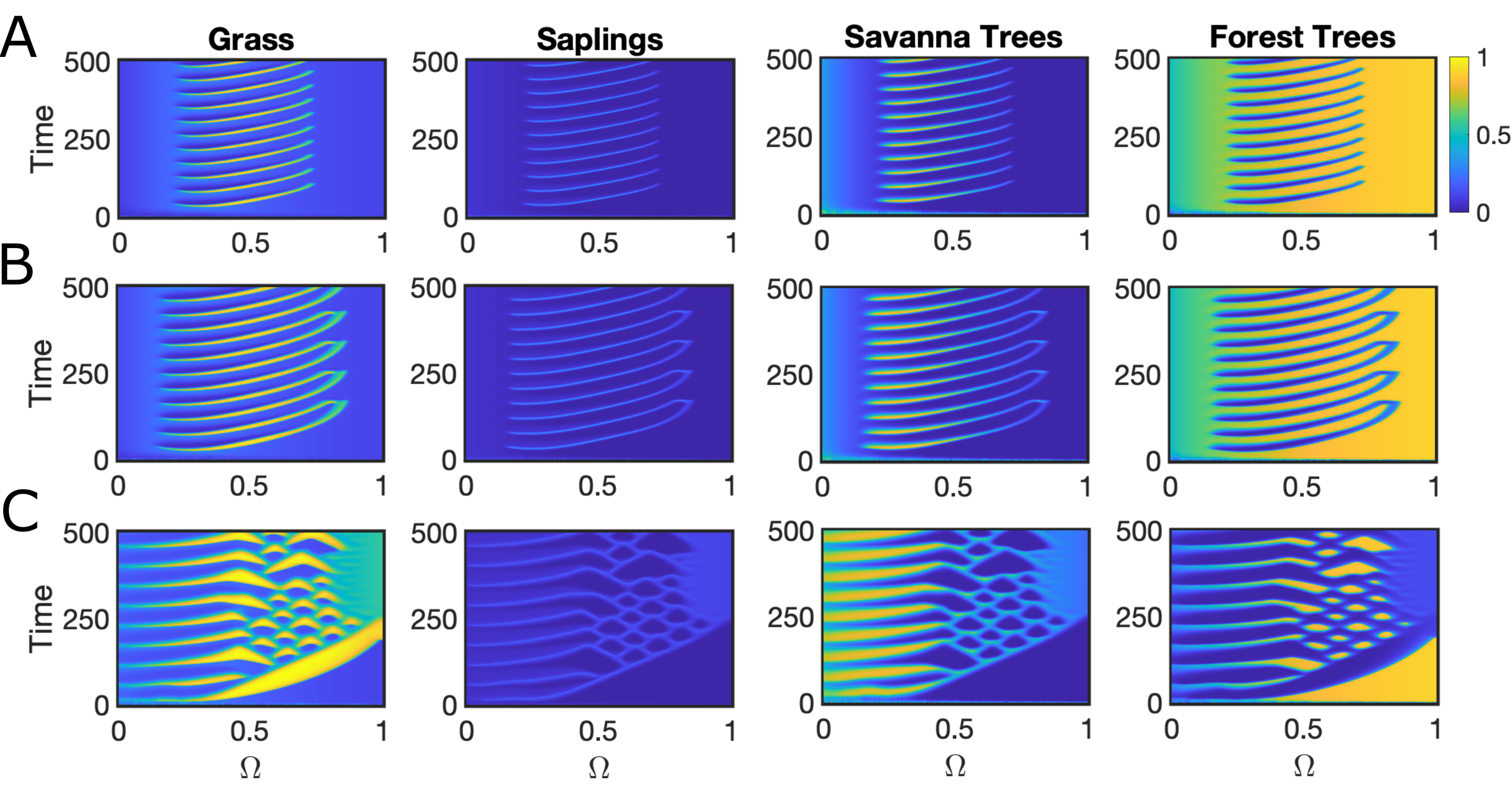}
	\caption{\textbf{A:} Simple (period 1) waves of invasion of savanna into forest led by a grass/fire front. \textbf{B:} More complex waves of invasion (period 2). \textbf{C:} Solution from the chaotic regime.}\label{fig.wave_types}
\end{figure}

Figure \ref{fig.wave_types} shows solutions for three different rainfall gradients which intersect the stable oscillations region of the nonspatial bifurcation diagram shown in Figure \ref{fig.codim2} B. The solution in Figure \ref{fig.wave_types} A varies from forest $\approx 0.5$ on the extreme left of the domain up to forest $\approx 1$ on the right but this forest dominance is punctuated by waves of invasion by the other \CS{functional types}. These waves are triggered by a wave of grass (carrying fire), which burns the forest trees and allows savanna to outcompete grass in the wake of the initial wave. The savanna  then becomes vulnerable to replacement by forest as it eliminates the grass\CS{, which carried fire to suppress the forest}. This cycle repeats in a similar manner to the simple periodic oscillations observed in the nonspatial Staver-Levin model~\cite{touboul2018complex}. The left and right endpoints of the gradient anchor the system in the savanna and forest states respectively, although the ``savanna state'' has unusually high forest tree presence. The solution shown in Figure \ref{fig.wave_types} B is a \SL{period-two} wave that appears to be caused when we enlarge the region of space which the gradient spends in the stable oscillation region of the nonspatial system. When this region is large enough, waves \SL{take} sufficiently long to reach the \SL{right-hand} termination point that another wave has already begun on the left side of the domain, eventually leading to \SL{period-doubling} due to the interaction of the two waves. Finally, in Figure \ref{fig.wave_types} C, we enlarge the oscillating region that the gradient passes through yet further and observe what appears to be a complicated quasi-periodic behavior or spatio-temporal chaos in the solution. The solutions in Figure \ref{fig.wave_types} were computed using reflecting boundary conditions but the solutions are qualitatively similar with open boundary conditions so the phenomena shown are not sensitive to this choice \edit{(see the Supplementary Materials for the corresponding simulations)}. 
\begin{figure}[h]
	\includegraphics[width=\textwidth]{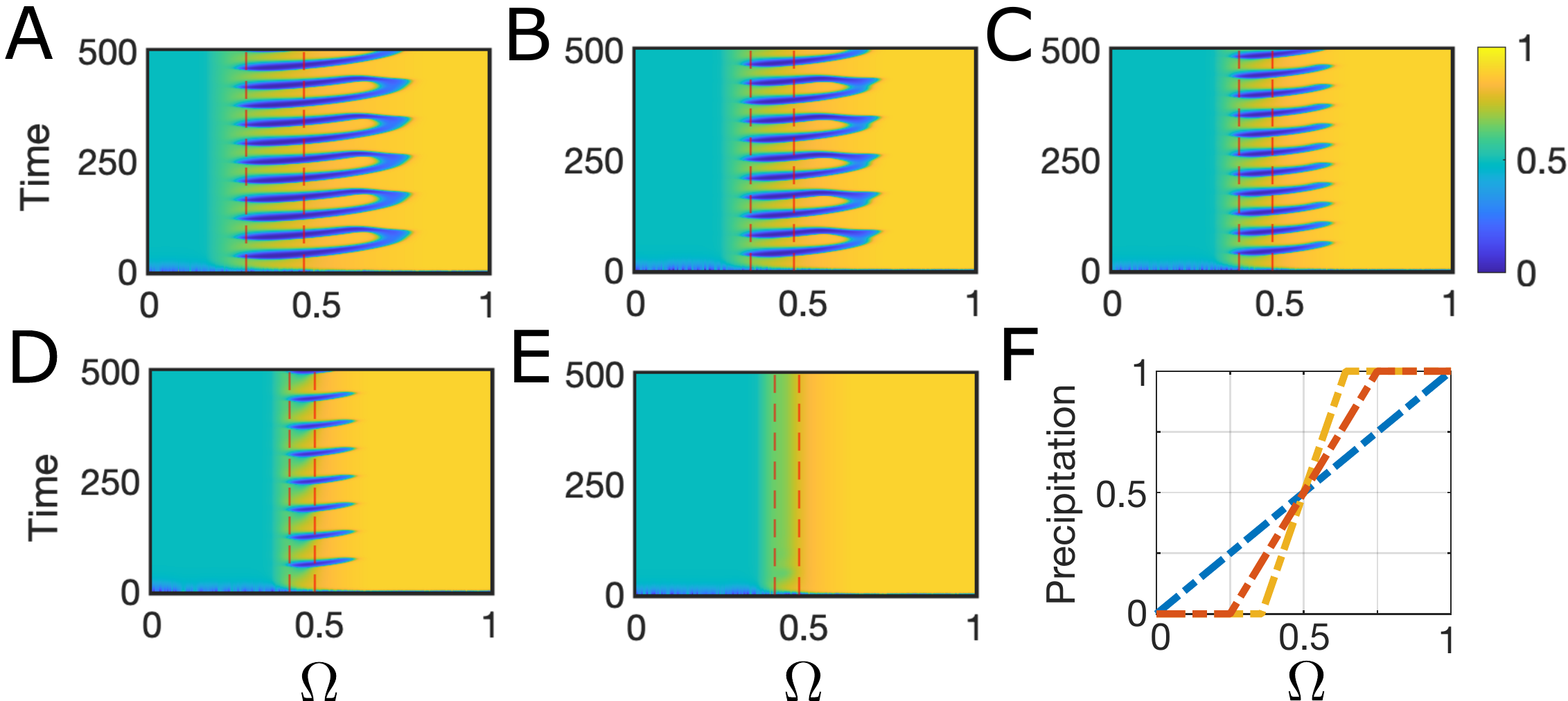}
	\caption{\edit{\textbf{A-E:} Forest component of each solution with the same gradient as Figure \ref{fig.wave_types}B but for varying slopes ($\sigma = 0.02$ in all cases). The slopes in panels A-E are: A - $1.5$, B - $2$, C - $2.5$, D - $3.4$, E - $3.42$.  \textbf{F:} Illustration of how the gradient slopes change across the other panels (blue gradient with slope $1$, red gradient with slope $2$, and yellow gradient with slope $3.5$)}.}\label{ref.waves_varying_slope}
\end{figure}
To further investigate the phenomenon of wave generation shown above, we considered the gradient from Figure \ref{fig.wave_types} B once more, but this time we allowed the slope of the gradient to vary in order to change the region of the spatial domain in the oscillatory regime (slope changes are illustrated in Figure \ref{ref.waves_varying_slope} F). Panels A, B and C of Figure \ref{ref.waves_varying_slope} show the forest component of the solution for slopes of $1$, $2$ and $2.5$ respectively. As the slope of the gradient increases, we soon revert from \SL{period-two} waves to simpler \SL{period-one} waves. The left vertical red dashed line in panels A-F of Figure \ref{ref.waves_varying_slope} indicates the onset of oscillations without spatial interactions, corresponding to the start of the fold of limit cycles associated with the subcritical Hopf curve marked by a solid dark green line in Figure \ref{fig.codim2} B; the right vertical red dashed line indicates the offest of oscillations, corresponding to the heteroclinic-to-saddle connection marked by a dashed light blue line in Figure \ref{fig.codim2} B. Even as the region of oscillations predicted by the nonspatial model shrinks, the waves continue to travel remarkably far past the right hand red line marking the non-spatial heteroclinic. Figures \ref{ref.waves_varying_slope} D and E show solutions for slope values of $3.4$ and $3.5$, with oscillations finally abating in panel E as the oscillatory region finally becomes too small to support stable oscillations in the spatial model. This example illustrates the phenomenon of transient passage through a family \SL{of} periodic orbits in a spatial model; it shows that real systems may pass through regimes \SL{that} support oscillations without spatial interactions but that these oscillations will only be seen in \SL{spatially-extended} models if the gradient passes through this region of parameter space sufficiently slowly. It is evidently of both theoretical and applied interest to understand in more detail when transient passage through oscillations will result in periodic solutions in the spatial system and how these periodic solutions may further bifurcate, as they do above in Figure \ref{fig.wave_types}. Immediate questions and challenges in this domain include identifying criteria to determine the onset or offset of spatio-temporal oscillations and characterizing transitions between different oscillatory regimes (e.g. period doubling, quasi-periodic and chaotic behaviors) for representative classes of spatial operators.

\section{Transient patterns at the frontiers: Applications to brain development}\label{sec.brain_patterning}
Physiological formation of cortical regions relies on the precise positioning of sharp and regular boundaries during embryonic development, \SL{which} are thought to be guided by the presence of positional cues (typically, gradients of morphogens)~\cite{flanagan:06,kiecker2005compartments}. Failures to forming cortical regions with sharp and regular boundaries at specified location in brain were implicated in many serious pathologies~\cite{hyman1991some,pederick2018abnormal,watrin2015causes}. 
Each territory is characterized by the expression of a specific combination of molecular marks that compete with each other and according to the gradients of morphogens promoting or repressing specific genes~\cite{oleary:07}. As an example, in the mammalian forebrain, cells arising from divisions of a common populations of progenitors are fated to become medial entorhinal cortex (MEC) cells or neocortical (NC) cells based on their exposure to patterning transcription factors (TFs) expressed in a graded fashion. In \cite{feng2021coup}, it was shown that changes in morphogen gradients not only led to simple boundary shifts, but also to the possibility of a shattered boundary with the formation of regular patterns of ectopic cortical domains. A mathematical model introduced in this paper proposed that this phenomenon could be related to a slow passage through a pattern-forming instability that we explore here in more detail.

\subsection{A model of brain arealization}
The model proposed in Feng et al. \cite{feng2021coup} was based principally on the following three experimentally and biologically motivated mechanisms:
\begin{enumerate}[(I.)]
	\item \emph{Competition:} The differentiation into MEC or NC neuronal identity is a competitive process (e.g.,competition on the genetic resources).
	\item \emph{External cues (heterogeneity):} Extracellular TF gradients favor differentiation into MEC cells on the posterior side and into NC cells on the anterior side of the brain.
	\item \emph{Aggregation/differential adhesion mechanisms:} Neurons, having a tendency to diffuse slowly, also have the ability to aggregate preferentially with cells of their own type (Feng et al. \cite{feng2021coup} showed that differential adhesion was the primary driver of cells preferentially aggregating with other cells of their own type).
\end{enumerate}
We thus arrive at the following continuum mathematical model of MEC and NC identity marker levels across the cortex:
\begin{subequations}\label{eq.KS2_sat}
	\begin{alignat}{2}
		\partial_t E(x,t) &= E(1 - E - k_1 N)  + D_E \Delta E - {\chi_1} \,\nabla \cdot \left({\Phi(E)} \,\nabla C_E\right) + \rho_E(x), \\
		\partial_t C_E(x,t) &= E - C_E + D_{C_E} \Delta C_E, \\
		\partial_t N(x,t) &= N(1 - N - k_2 E)  + D_N \Delta N - \chi_2 \nabla \cdot \left( \Phi(N)\nabla C_N\right) + \rho_N(x) , \\
		\partial_t C_N(x,t) &= N - C_N + D_{C_N} \Delta C_N,\quad (x,t)\in\Omega\times \mathbb{R}^+,
	\end{alignat}
\end{subequations}
where the species $E$ accounts for MEC cell fate markers, $N$ accounts for NC cell fate markers, $\nabla$ denotes the gradient operator, $\nabla \cdot $ denotes the divergence operator and $\Delta$ is the diffusion operator. The competition between $E$ and $N$ is reflected by classical Lotka-Volterra competitive reaction dynamics; we choose $k_1=k_2=2$ to inhabit a bistable regime mirroring the common lineage of MEC and NC progenitor cells when spatial interactions are neglected. External signals promote the expression of each gene: $E$ is promoted at a rate $\rho_E(x)$ at position $x$ (COUP-TFI effects and other signals promoting MEC fate) and $N$ is promoted at a rate $\rho_N(x)$ (TF effects promoting NC fate) at position $x$. We account for cell aggregation via two monitor species, $C_E$ and $C_N$, generated respectively by cells expressing $E$ or $N$ markers and subject to degradation and diffusion. $C_E$ and $C_N$ attract $E$ and $N$ cells respectively with $\chi_1$ and $\chi_2$ denoting the strengths of the aggregation forces. Considering our intended application, it makes sense to always consider \eqref{eq.KS2_sat} with no-flux boundary condition for $\Omega\subset\mathbb{R}^n$ with $n \leq 3$. It should be stressed that our model is phenomenological in nature and we have not aimed at a detailed biophysical description of the physical and cellular process involved, but rather we have tried to capture their effects qualitatively in a parsimonious and interpretable model.

The model \eqref{eq.KS2_sat} follows the classical framework of Keller-Segel chemotaxis systems with cell aggregation limited by saturation effects at high cell density by $\Phi:\mathbb{R}\mapsto \mathbb{R}$~\cite{hillen2009user,keller1970initiation,keller1971model}. Mathematically, this makes the model more stable and avoids the potential for finite-time blow-up of solutions, as often arises in other Keller-Segel-like models (cf. \cite{hillen2001global}). We employ a Ricker-type saturation function, i.e. 
$\Phi(z)= \alpha z e^{-\alpha z}$ with $\alpha>0$, but other choices appropriate for capturing saturation produce qualitatively similar results to that shown below and in Feng et al.~\cite{feng2021coup}. 

\edit{Table~\ref{table.param} summarizes the parameters, their biological interpretations, and their default values in the brain patterning model.}
\begin{table}[h]
		\centering
	\caption{Summary of parameters for the brain arealization model}
	\centering
	\edit{
	\begin{tabular}{@{}lcc@{}}
		\toprule
		Biological interpretation    & Expression & Default Value \\
		\midrule
		Entorhinal marker level & $E $    &  [dynamic] \\
		Entorhinal monitor concentration    & $C_E$   & [dynamic]    \\
		Neocortex marker level         & $N$   & [dynamic]    \\
		Neocortex monitor concentration     & $C_N$ & [dynamic] \\
		COUP-TF\RN{1} effects & $\rho_E(x)$ & [piecewise linear] \\   
		Patterning TF effects promoting neocortex & $\rho_N(x)$ & [piecewise linear]  \\   
		Entorhinal competition & $k_1$ & 2 \\   
		Neocortex competition & $k_2$ & 2 \\ 
		Aggregation saturation parameter  & $\alpha$ & 1.2 \\
		Entorhinal adhesion strength  & $\chi_1$ & 1.5 \\
		Neocortex adhesion strength  & $\chi_2$ & 1.5  \\
		Entorhinal cell diffusion & $D_E$ & 0.2 \\
		Entorhinal monitor diffusion & $D_{C_E}$ & 0.2 \\
		Neocortex cell diffusion & $D_N$ & 0.2 \\
		Neocortex monitor diffusion & $D_{C_N}$ & 0.2 \\
		\bottomrule
	\end{tabular}\label{table.param}
}
\end{table}

\subsection{Transient passage through a pattern forming instability}

Our primary goal in this paper is to illustrate the mathematical mechanism by which the model \eqref{eq.KS2_sat} explains the experimental observations of ectopic MEC from Feng et al.~\cite{feng2021coup} and hence we refer the interested reader to that paper for further biological details. Indeed, with appropriately chosen TF gradients (i.e. $\rho_E(x)$ and $\rho_N(x)$), the system \eqref{eq.KS2_sat} qualitatively matches all observed brain patterning phenotypes from the Feng experiments~\cite{feng2021coup}. The only modification we make to the model \eqref{eq.KS2_sat} is to choose $\chi_1$ and $\chi_2$ as fixed constants. It turns out that these coefficients are also influenced somewhat by COUP-TFI and this dual role of COUP-TFI is at the heart of the complex boundary shattering phenomenon observed in vivo, but this is not necessary to explain the dynamics of interest here from a mathematical perspective.  

\begin{figure}[h]
	\centering
	\includegraphics[width=0.7\textwidth]{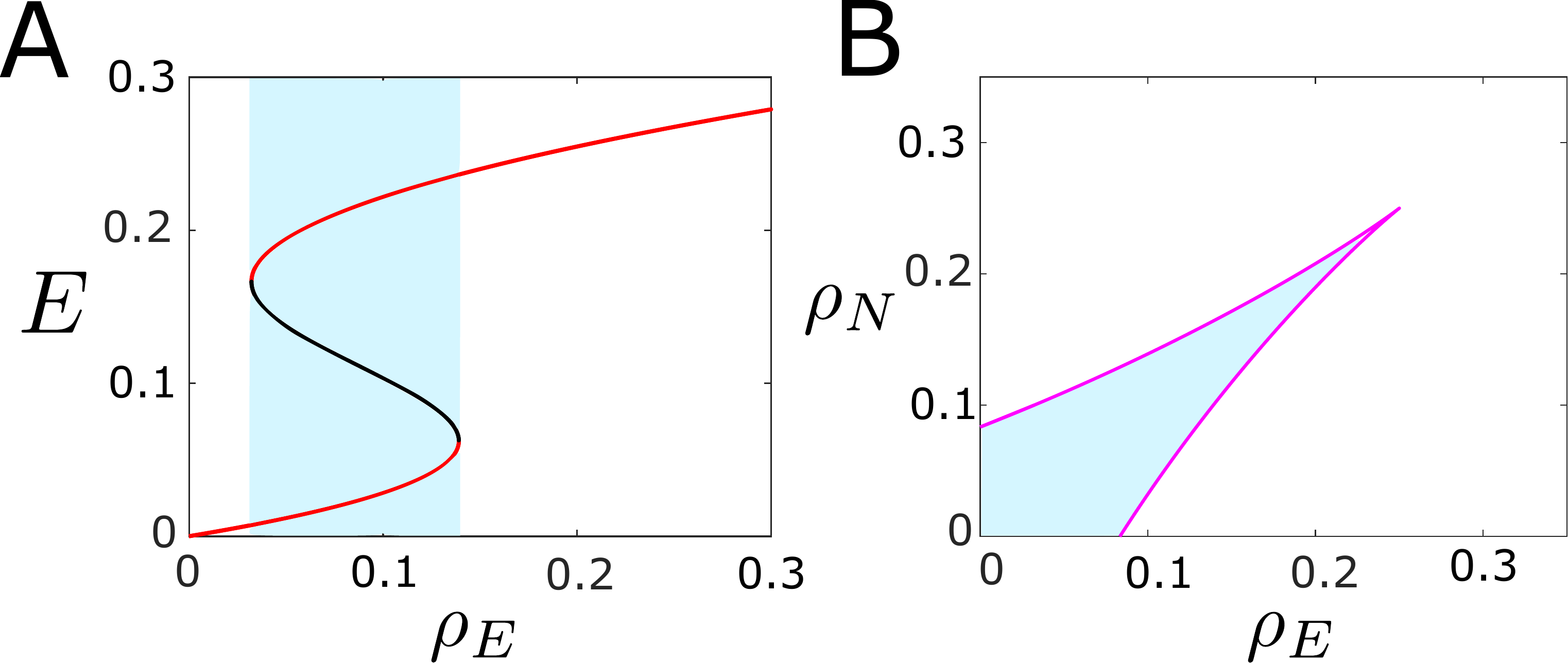}
	\caption{\edit{\textbf{A:} One-parameter bifurcation diagram for solutions of \eqref{eq.KS2_sat} without spatial interactions with $\rho_N = 0.1$. \textbf{B:} Two-parameter bifurcation diagram for \eqref{eq.KS2_sat} without spatial interactions. Magenta lines are curves of saddle-nodes which collide in a cusp bifurcation.}}\label{fig.nonspatial_brain}
\end{figure}

First, consider the system \eqref{eq.KS2_sat} posed on a homogenous domain, i.e. $\rho_E(x) \equiv \rho_E \in\mathbb{R}^+$ and $\rho_N(x) \equiv \rho_N \in\mathbb{R}^+$. Homogeneous solutions $(\bar{E},\bar{C}_E \bar{N},\bar{C}_N)$ to \eqref{eq.KS2_sat} obey the following system of nonlinear equations:
\begin{subequations}\label{eq.KS2_sat_equilibria_2}
	\begin{alignat}{2}
		0 &= \bar{E}(1 - \bar{E} - 2 \bar{N}) + \rho_E, \label{eq.3a}\\ 
		0 &= \bar{N}(1 - \bar{N} - 2 \bar{E})  + \rho_N,\label{eq.equilibrium2}
	\end{alignat}
\end{subequations}
with $\bar{C}_E = \bar{E}$ and $\bar{N} = \bar{C}_N$. The system \eqref{eq.KS2_sat_equilibria_2} has between one and three solutions depending on the values of $\rho_E$ and $\rho_N$. In the absence of spatial interactions, the stability of these equilbria is shown as a function of $\rho_E$ and $\rho_N$ in Figure \ref{fig.nonspatial_brain} panels A and B. Figure \ref{fig.nonspatial_brain} A tracks the stability of homogeneous solutions as a function of $\rho_E$ with $\rho_N=0.1$ and shows that there is a region of bistability (in light blue) with two stable homogeneous solutions separated by an unstable equilibrium. In Figure \ref{fig.nonspatial_brain} B we allow both $\rho_E$ and $\rho_N$ to vary and we observe that the two saddle-node bifurcations from panel A collide in a codimension 2 cusp bifurcation at $(0.25,0.25)$ in $\rho_E$-$\rho_N$ space. Between these \SL{saddle-node} curves there is a significant region of bistability (in light blue) while outside of this region there is a single solution \SL{that} is $N$ dominated above the line $\rho_E=\rho_N$ and $E$ dominated below that line.

\begin{figure}[h]
	\centering
	\includegraphics[width=\textwidth]{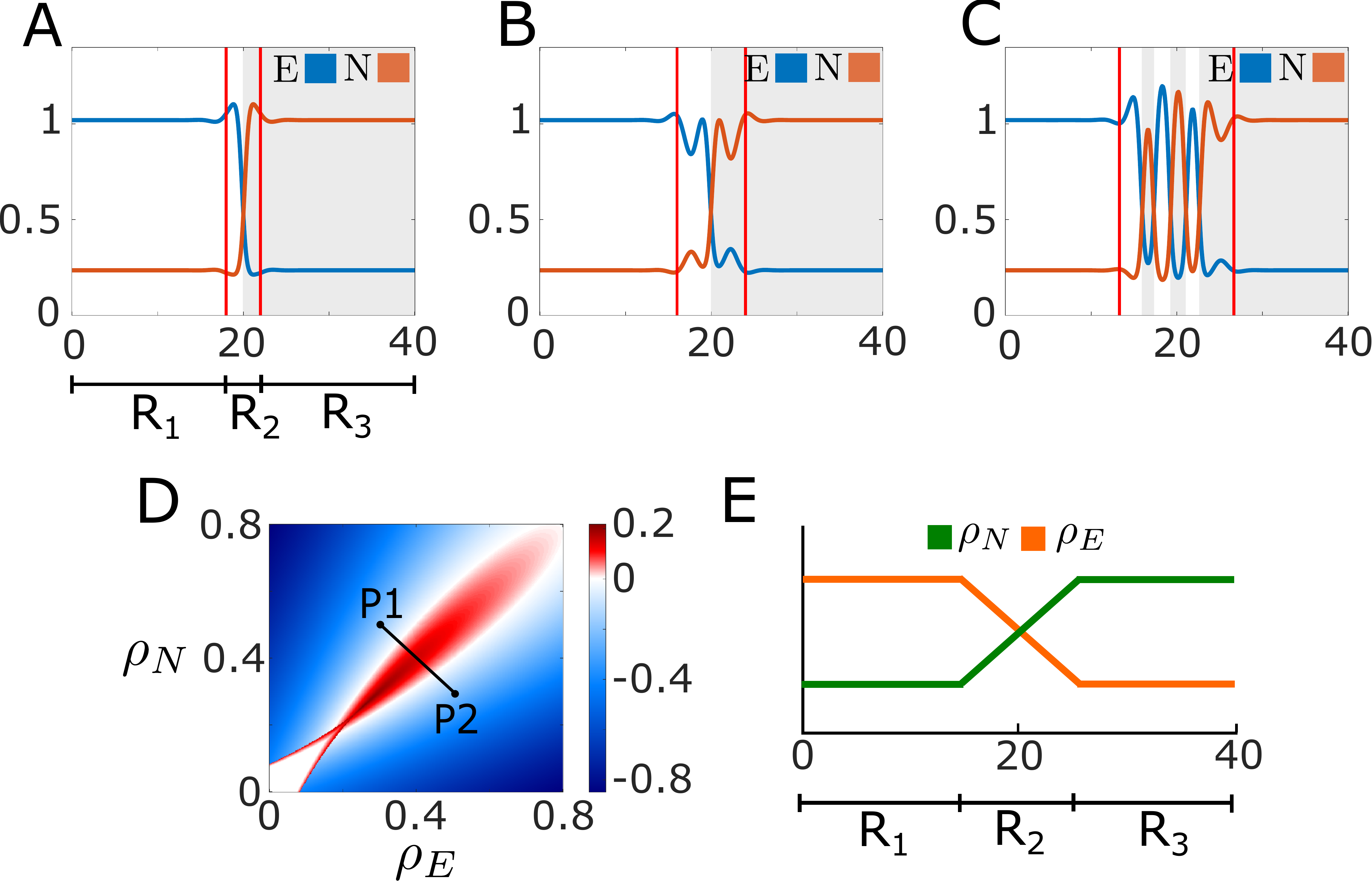}
	\caption{\edit{\textbf{A-C:} Model simulations with a linear morphogen gradient in the region marked by red vertical lines; constant morphogen levels fixed at the gradient endpoints outside the red marked region. All solutions shown at time $t = 200$. \textbf{D:} Heatmap of the maximum principal eigenvalue obtained when linearizing around all homogeneous equilibria that are stable without spatial interactions (diffusion coefficients all $0.2$ and $\chi_1=\chi_2=1.5$). \textbf{E:} Illustration of how we vary the gradient in $\rho_E$-$\rho_N$-space across the 1D spatial domain $\Omega=[0,40]$.}}\label{fig.1D_linear}
\end{figure}

Turning now to the full spatial model, we can linearize about the homogeneous equilibria discussed above and classical linear stability analysis reveals pattern forming instabilities will emerge for a range of values of the diffusion coefficients and chemotactic strength parameters (see Feng et al. \cite{feng2021coup} supplementary information for details). Figure \ref{fig.1D_linear} D shows the potential instability region in $\rho_E$-$\rho_N$-space via a heatmap of the maximum of the principal eigenvalues of the linearized operators about each \SL{of} the homogeneous equilibria \SL{that} were stable in the absence of spatial interactions; this red region is typically referred to as the Turing space~\cite{murray1982parameter}. Numerical simulations of solutions to \eqref{eq.KS2_sat} confirm the presence of patterns in the red regions of Figure \ref{fig.1D_linear} D, and reveal the nature and diversity of these patterns (spots, mixed spot-stripe and labyrinths patterns were observed for this system in 2D domains). 

The analysis above is valid for $\rho_E$ and $\rho_N$ equal to some constant fixed values but, biologically, they both vary across the cortex. In our idealized model, we assume $\rho_E$ simply promotes $E$ and $\rho_N$ promotes $N$. Moreover, $E$ should dominate at one side of the spatial domain, while $N$ dominates at the other. Thus a reasonable choice of the gradient must start above the line $\rho_E = \rho_N$ (which represents perfectly balanced competition) and end below it. One such choice is shown in Figure \ref{fig.1D_linear} D by the solid black line giving a linear gradient varying from points $P1$ to $P2$ and cutting through the red instability region. However, as alluded to earlier in the discussion of spatio-temporal oscillations in section \ref{sec.rainfall_gradient}, the key question in terms of the resulting spatial dynamics is: how quickly does this gradient cross the pattern forming region? In Figure \ref{fig.1D_linear} E, we illustrate how we can vary the slope of the gradient or in other words, vary the length of the region $R_2$ to adjust how long the gradient spends in the pattern forming regime from the homogeneous domain problem. Figure \ref{fig.1D_linear} panels A, B and C show three solutions of the model \eqref{eq.KS2_sat} using the piecewise linear gradient from Figure \ref{fig.1D_linear} D on the 1D spatial domain $\Omega=[0,40]$ with the edges of the region $R_2$ marked by red vertical lines in each case. In Figure \ref{fig.1D_linear} A, we don't observe standard pattern formation but rather a type of front pinned solution similar to those discussed earlier for the SL model (see section \ref{sec.rainfall_gradient}). In panel B, we see the solution form a couple of small amplitude spikes, coming closer to standard pattern forming behavior, and in panel C, the solution now supports a multiple spikes and resembles the expected solution for a pattern forming system in the subdomain $R_2$. 

Notably, the spikes in Figure \ref{fig.1D_linear} F do not have a single wavelength or frequency. Multiple wavelengths appear to be present as we cross the instability region; this is distinct from a standard Turing pattern with a single dominant wavelength, and analogous to the model and experiments reported in~\cite{feng2021coup} where there appears to be a blend of spot, stripe and labyrinthine patterns, mirroring exactly the phenomenology observed in vivo. Clearly, the width of $R_2$ crucially influences the ability of the system to form patterns and the nature of those patterns. Moreover, the wavelength of the patterns supported by the system posed on a homogeneous domain relative to the width of $R_2$ appears important in determining whether patterns can ``fit'' into potential patterning region. 

To the best of our knowledge, the biologically motivated dynamics above, which we term ``transient passage through a pattern forming instability'', do not fit into any of the established frameworks that provide criteria for the onset of Turing-like patterns, not even those designed to incorporate heterogeneity (e.g. \cite{kozak2019pattern,krause2020one,van2021pattern}). Moreover, it would evidently be of interest to characterize the nature of the patterns further to understand if different qualitative patterns are possible (for two or three dimensional domains), as these may have different implications for applications in biology. Thus, this applied example provides strong motivation for the continued development of mathematical tools for studying pattern formation in the presence of spatial heterogeneities, particularly for non-reaction-diffusion systems and large continuously-varying heterogeneities.

\section{Conclusions and Discussion}
\edit{Despite the prevalence of pattern formation in heterogeneous domains in nature, the abundant theoretical literature in pattern formation has still not developed the tools to address these questions for general transient passages through instabilities. The two biologically motivated examples introduced here do not fit into any of the established frameworks that provide criteria for the onset of Turing-like patterns, not even those designed to incorporate heterogeneity. We hope that highlighting these examples will motivate new mathematical theory to better understand the underlying dynamics in each case.}

\edit{Existing work to allow heterogeneity in pattern forming systems has achieved great progress in characterizing small heterogeneous perturbations of homogeneous systems~\cite{benson1998unravelling} or bifurcations of heterogeneous systems with steady states that are explicitly derived. There has also been significant work in the neuroscience literature studying more complex heterogeneities as inputs to neural field equations~\cite{faye2014pulsatile,kilpatrick2008traveling,kilpatrick2013optimizing,kilpatrick2013wandering}, and some recent work in theoretical ecology highlighting the complexity of dynamics that non-monotonic gradients can generate in bistable systems~\cite{bastiaansen2022fragmented}. In contrast, we are concerned with characterizing the properties of solutions to spatially heterogeneous systems as they relate to a bifurcating behavior in the underlying homogeneous (typically non-spatial) system. In other words, the mechanism of interest here is rather spatial models with multistable behavior which have spatial gradients (heterogeneities) connecting regions with different underlying dynamical behaviour. Thus, the complexity in the dynamics emerges from the rate at which (monotonic) gradients connect regions with different dynamic behaviors. For example, a much studied paradigm for the two species SL model (forest-grass) presented in Section \ref{sec.rainfall_GF} is one in which there are essentially three regions of the spatial domain: a region where the grassland state is stable, a region of bistability between forest and grass, and a region where forest is stable. Similarly, in the brain arealization model, there is a spatial gradient crossing from a region of enthorinal fate dominance to a region of neocortical fate dominance. However, in this case, the complex dynamics emerge when we increase the adhesion strength sufficiently to introduce a third distinct region in which pattern formation occurs, leading to what we term \emph{transient passage through a pattern forming instability}.}

\edit{Mathematically, we expect some asymptotic regimes to be amenable to analysis. In particular, we expect that situations with very localized interactions (e.g., vanishing diffusion for instance) will closely match the non-spatial dynamics away from bifurcations and in regions of space with unique homogeneous stable attractors, with rapid transitions between distinct patterns at bifurcations or within multi-stable regions. Existing mathematical works in this domain have focused on competition models that may include bistability, and only recently qualitative properties in asymptotic regimes were derived~\cite{perthame2015competition}. These situations are in fact only among the simplest cases of the general question of pattern formation in heterogeneous domains and through bifurcations. Mathematically, we expect that tools from perturbative analysis, geometric (singular) perturbation theory~\cite{bastiaansen2020pulse}, or changes of parameterization of solutions as used in~\cite{perthame2015competition} to derive viscosity solutions could allow rigorous characterizations of those regimes. Regimes with fixed diffusion but very slow variations of the environment in space are likely similar, and will sometimes be precisely matching regimes of very localized interactions in space through appropriate changes of variables. Another asymptotic regime that will likely be amenable to analysis is systems with sigmoidal gradients of heterogeneity in the limit of very sharp gradients. In these regimes, we expect to observe a convergence towards the solution of a system with patchy heterogeneity  (Heaviside step function) with a single transition between the leftmost and rightmost regimes and no impact of the transient patterns, as in Figures \ref{ref.waves_varying_slope}E and \ref{fig.1D_linear}A, and as studied in~\cite{kozak2019pattern,page2003pattern}. For all these questions, systems with multiple spatially homogeneous solutions (and possibly transient equilibria in space) will constitute the first models to analyze.}

Systems with spatially transient oscillations or patterns will be associated with a richer phenomenology. In systems with spatially transient Turing-like patterns, we expect not only the diffusion properties to play a role in the emergence of a pattern at the transition, but also the intrinsic properties of the pattern itself (as its length scale or the modes of instability associated). Indeed, we expect that contrasting with transient spatially homogeneous regimes, the emergence of a pattern at the transition will also depend on whether or not multiple patterns ``fit'' in the region of instability, and possibly not fully expressed patterns may emerge as well as in Figure \ref{fig.1D_linear} A. In multiple dimensions, the nature of the pattern and progressive morphing of it along the gradient will also arise, and this already plays a role in the patterns observed both from numerical simulations of the model and the corresponding experiments in \cite{feng2021coup}. In systems with spatially transient oscillations, not only the question of the mere existence of an oscillation at the transition arises, but problems related to the waves generated, their regularities and their bifurcations have a strong impact on the qualitative behavior of the solution. The prospects for theoretical results in this case seems somewhat more pessimistic for nonlocal spatial models, such as the SL model studied here, but it is likely that this phenomenon can fruitfully be studied for reduced form phase models with simple coupling structures. Moreover, while the theoretical predictions related to the wave-like solutions of the SL model presented here are not easily testable, similar mechanisms could be tested in smaller scale experimental systems, such as chemical oscillators~\cite{crowley1989experimental}.

\appendix

\section{The Staver-Levin Model}
\subsection{Parameters and numerics}\label{sec.numerical_parameters}
All codes to generate the figures from this paper are maintained on Github at \\ \href{https://github.com/patterd2/SL_model_rainfall_gradient}{github.com/patterd2/SL\_model\_rainfall\_gradient}.

The parameter gradients used in the paper for $\alpha$ and $\beta$ are linear, i.e.
\begin{equation}
	\alpha(x) = \alpha_c + \alpha_s x, \quad \beta(x) = \beta_c + \beta_s x, \quad x \in \Omega = [0,1].
\end{equation}
\begin{table}[H]
	\centering
	\begin{tabular}{lcccccccccc}
		Parameter & $\alpha_c$ & $\alpha_s$ & $\beta_c$ & $\beta_s$   \\ \hline
		Figures 1 \& 2 & 0.5 & 1.25   & N/A        & N/A    \\
		Figure 4 & 0.8 & 0.5   & 0.15        & 0.1 \\
		Figure 5A & 0.2 & 0.8   & 1.9        & 0.1    \\
		Figures 5B \& 6 & 0.2 & 0.8   & 1.5        & 0.1   \\
		Figure 5C & 0.2 & 0.8   & 0.4        & 0.1
	\end{tabular}
	\caption{Parameter values for the gradients for each figure from the main text.}
	\label{table.gradients}
\end{table}
In Figure 6 of the main text we adjust the speed at which we move along the rainfall gradient in the following way. We choose:
\begin{equation}
	\alpha(x) = \alpha_c + \alpha_s P(x), \quad \beta(x) = \beta_c + \beta_s P(x), \quad x \in \Omega = [0,1].
\end{equation}
where for slope parameter $P_s \geq 1$, we define
\[
P(x) = 
\begin{cases}
	0, \quad &  x \in \left[0,\,(1-1/P_s)/2\right), \\
	P_s \, (x - 0.5), \quad &  x \in \left((1-1/P_s)/2 ,\, (1+1/P_s)/2 \right), \\
	1, \quad &  x \in \left((1+1/P_s)/2,\,1\right].
\end{cases}
\]
We refer to $P_s$ as the ``slope parameter'' in the main text. The solution shown in Figure 6 A has $P_s = 1.5$, panel B is for $P_s = 2$, panel C is for $P_s = 2.5$, panel D is for $P_s = 3.4$ and panel E is with $P_s=3.42$.

We computed approximate numerical solutions to the spatial two-species SL model, given by equation \eqref{eq.GrassForest_space}, by first discretizing time with an explicit Euler scheme to obtain
\begin{equation}
	\begin{split}\label{eq.Appendix_euler}
		G(x,n+1) = G(x,n) +  h \left\{ \left(1 - G(x,n)\right) \,\phi\left( \int_{\mathbb{R}} w(x-y)\,G(y,n)\, dy \right) \right. \\
		\left. - G(x,n) \,\left(1 - \int_{\mathbb{R}} J_F(x-y)\, G(y,n) \,dy \right) \vphantom{\left(1 - G(x,n)\right) \,\phi\left( \int_{\mathbb{R}} w(x-y)\,G(y,n)\, dy \right)} \right\},\quad h>0,\quad (x, n)\in[0,1] \times  \mathbb{Z}^+.
	\end{split}
\end{equation}
We then discretize the integrals in \eqref{eq.Appendix_euler} using the 1D trapezoidal rule and approximate the solution on the evenly spaced grid $\{0,\Delta,2\Delta,\dots,5\}$ for some $\Delta>0$. In practice, we found that a step size $h<0.1$ and $200$ spatial grid points was sufficient to ensure numerical stability of the scheme and consistency of the qualitative dynamics, and the scheme remained stable as we decreased the time step and increased the number of grid points. Solutions shown in the main text are for 400 grid points. A similar discretization was used for the four-species SL model.

\subsection{Boundary Conditions for nonlocal operators}\label{sec.boundary_conditions}
To fix ideas, consider the (ill-posed) nonlocal operator given by
\begin{equation}\label{eq.boundary_operator}
	\mathcal{K}[u](x,t) = \int_\mathbb{R} J(x-y)u(y,t)\,dy,\quad t\in\mathbb{R}^+,
\end{equation}
where the function $u$ is defined on $[0,L]\times\mathbb{R}^+$. From an applied perspective, the value of $\mathcal{K}[u](x,t) $ represents the cumulative impact of seed dispersal (or fire signal) on $u(x,t)$, the vegetation present at position $x$ at time $t$. The value of the integrand in this formula, $J(x-y)u(y,t)$, is dispersal (either seed or fire transmission) from vegetation present at position $y \in \mathbb{R}$ but for numerical simulations the solution is only defined on a finite domain, i.e. $u(y,t)$ is only defined for $y \in [0,L]$ for each $t\in\mathbb{R}^+$. Hence we must ask: what choice of boundary conditions extends $u$ from a function on $[0,L]$ to a function on $\mathbb{R}$ and respects the character of the ecological application at hand?

\begin{figure}[h]
	\centering
	\includegraphics[width=0.7\linewidth]{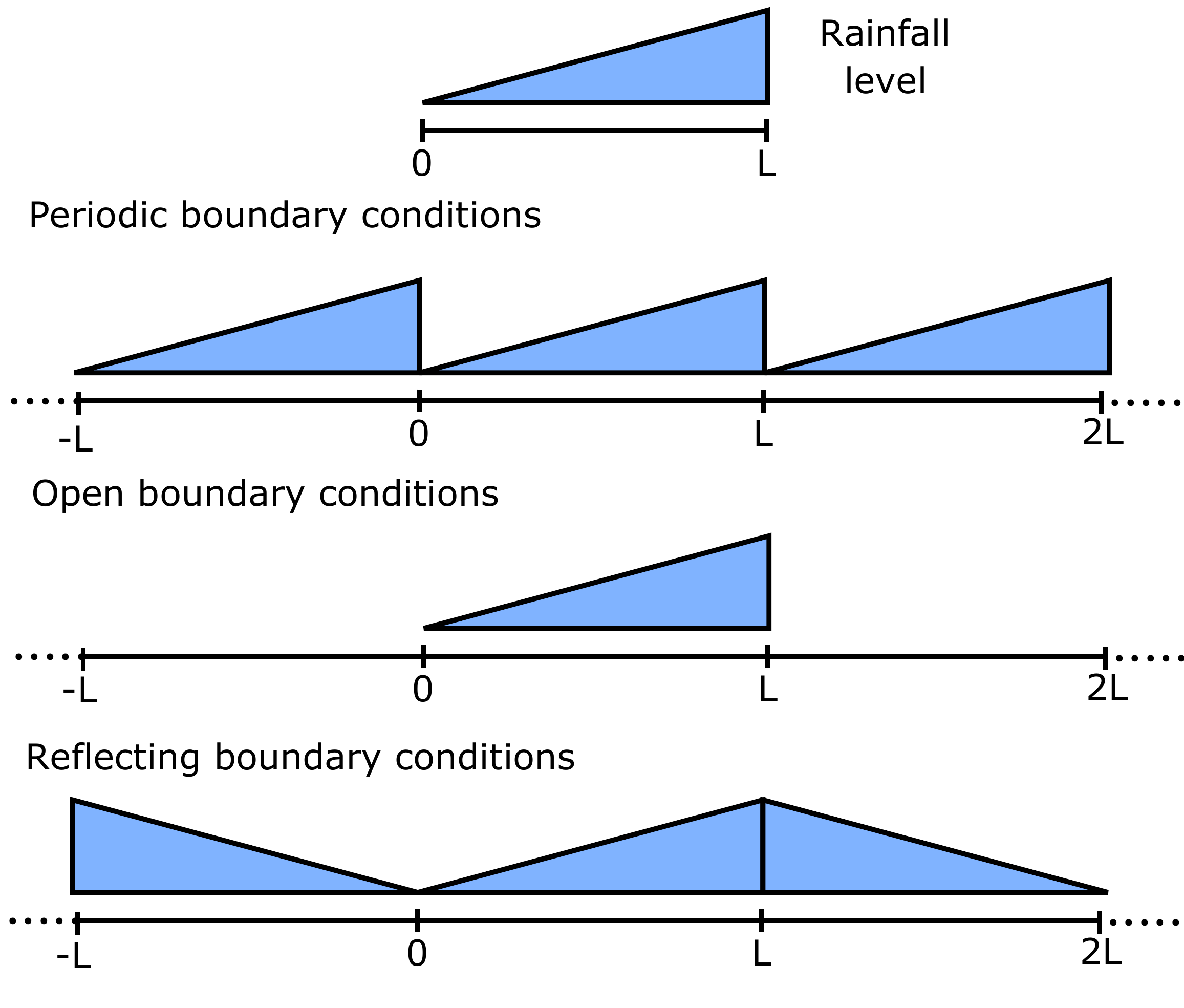}
	\caption{Comparison of the practical effects of the different boundary conditions on the heterogeneous medium structure in the rainfall gradient model.}\label{fig.boundary}
\end{figure}

\emph{Periodic boundary conditions:} Nonlocal operators similar in structure to \eqref{SL_integral} feature extensively in the neuroscience literature in so-called neural mass or neural field models, virtually always with \emph{periodic boundary conditions}. With periodic boundary conditions one takes the standard periodic extension of $u$ given by
\[
\tilde{u}(x,t) = u(x \text{ mod }L,t), \quad x \in \mathbb{R}.
\]
This type of boundary is typically favored in a spatially homogeneous problem in which boundary effects on the dynamics are not expected be significant. Although we do not expect significant boundary effects in our problem, periodic boundaries will not be a appropriate for the heterogeneous medium problem we wish to consider; it will also effectively make the heterogeneity itself spatially periodic and this can lead to significant undesirable boundary effect (e.g. forest at high rainfall invading into savanna at low rainfall via the boundary, cf. Figure \ref{fig.boundary}). 

\emph{Open boundary conditions:} The most intuitively appealing option is to take inspiration from the real-world and assume that seeds (resp. fire transmission) which exit the domain boundaries at $0$ and $L$ are simply removed from the system, i.e. the boundaries are ``open''. Additionally, no seeds or fires are transmitted into the domain from outside. We refer to this as an \emph{open  boundary condition} and it implies that we extend the definition of $u$ as follows:
\[
\tilde{u}(x,t) = \begin{cases}
	u(x,t), \quad &x\in[0,L], \\
	0, &x\notin[0,L].
\end{cases}
\]
This type of boundary condition has the advantage of being physically realistic for both seed dispersal and fire transmission but will introduce noticeable boundary effects in the solutions. If the dispersal kernels are sufficiently localized relative to the domain, these boundary effects will not qualitatively impact which stable solutions are selected. For sufficiently long range dispersal, stability of solutions can be effected by boundary effects with this boundary condition but this will not be relevant at the length scales we consider for our applications.

\emph{Reflecting boundary conditions}: This is a particular type of periodic extension of the solution which we employ in some of our numerical investigations because it minimizes boundary effects for even (symmetric) kernels. First extend $u$ from a function on $[0,L]$ to $[-L,L]$ by reflection:
\[
u_{R}(x,t) = \begin{cases}
	u(x,t), \quad &x\in[0,L],\\
	u(-x,t), &x\in[-L,0].
\end{cases}
\]
Now let $\tilde{u}$ be the standard $2L$ periodic extension of $u_R$ (as defined above) so that $\tilde{u}$ is defined on all of $\mathbb{R}$. The reflection symmetry introduced in this extension, plus the choice of even kernels, ensures that the solutions to the IDE system \eqref{SL_integral} closely approximate those of the underlying nonspatial model in the limit as dispersal tends to zero. This choice also has the virtue of putting ``wetter'' regions next to ``wetter'' regions in the context of the rainfall gradient model (see Figure \ref{fig.boundary}) and hence we choose to show solutions using this boundary condition in the main text.

\subsection{Non-monotonic gradients and robustness in the SL model}

\begin{figure}[H]
	\centering
	\includegraphics[width=0.8\textwidth]{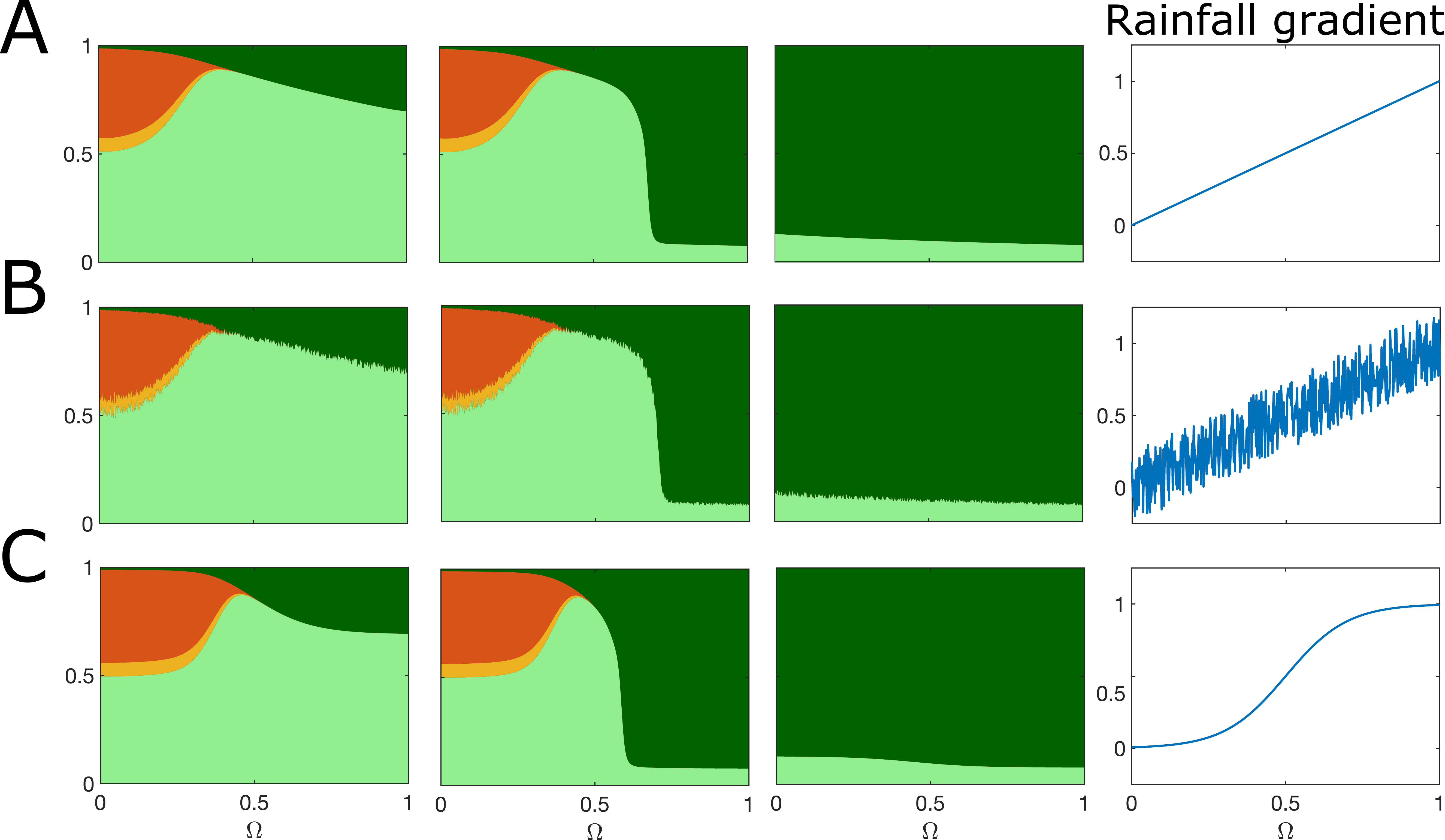}
	\caption{Row A shows the solutions from Figure 4 along with the corresponding linear rainfall gradient. Row B shows solutions for the same parameters but with noise added to the rainfall gradient. Row C shows solutions for a similar parameter regime but with a nonlinear, but still monotonic, rainfall gradient that starts and ends at the same points in $\alpha$-$\beta$ space.}\label{fig.nonmono_sup}
\end{figure}

In Figure \ref{fig.nonmono_sup} we show simulations corresponding to the setup of Figure 4 of the main text, but with non-monotonic rainfall gradients to emphasize the robustness of the multistability observed in this example. In addition to our default linear gradient, we also show qualitatively similar solutions for a (stochastic) noisy version of the linear gradient (row B) and a nonlinear gradient with the same start and end points (row C). These simulations were carried out with reflecting boundary conditions.

\section{Wave-like solutions to the SL model under different boundary conditions}\label{sec.waves_supp}
Figure \ref{fig.waves_sup} below shows more detail on the potentially chaotic solutions observed in Figure 6 of the main text. In particular, the spatial averages of the solution components show no periodic character, even over a very long time interval.
\begin{figure}[H]
	\centering
	\includegraphics[width=0.9\textwidth]{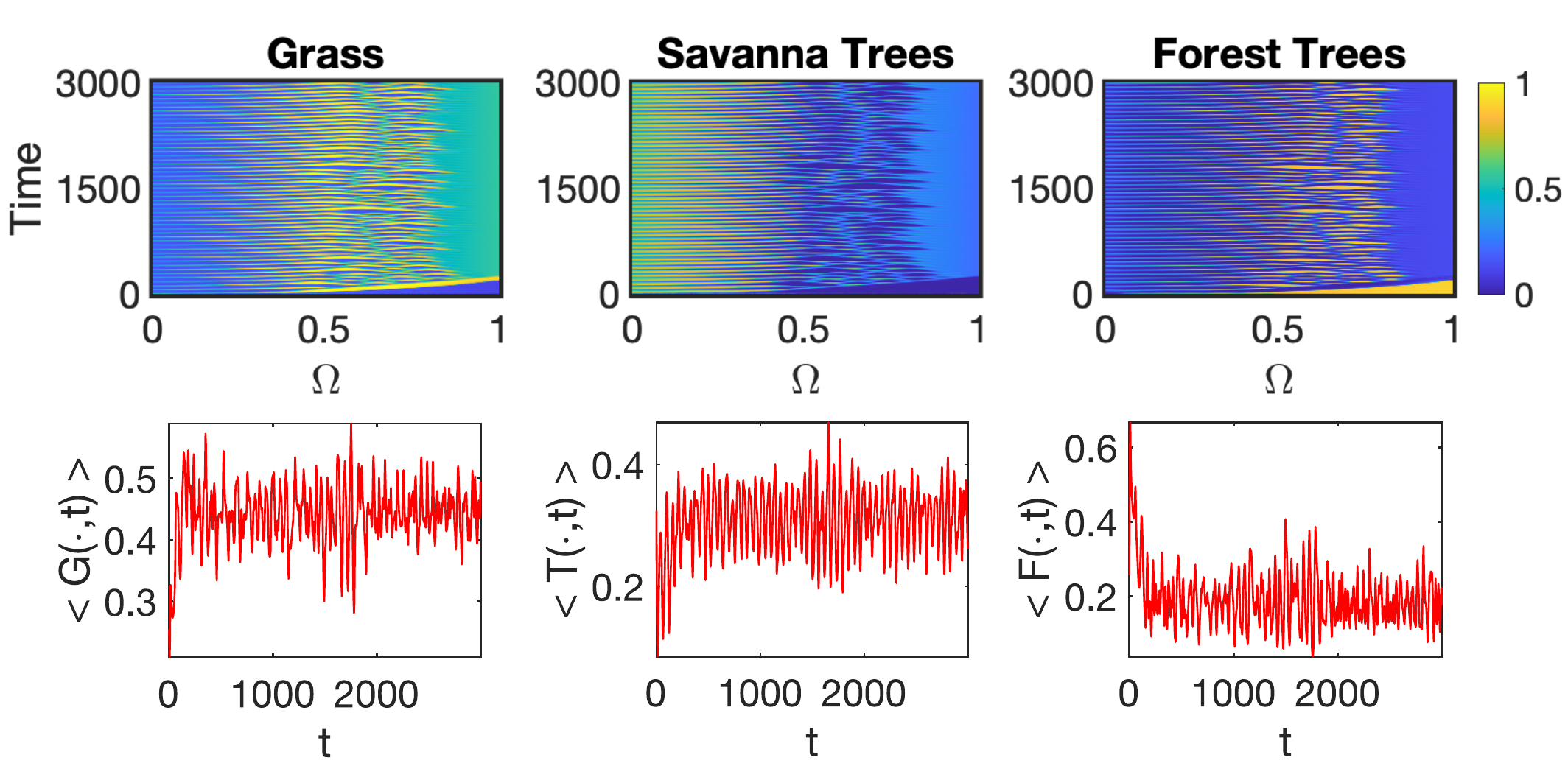}
	\caption{Row 1 shows the solutions from Figure 6C over a longer time horizon. Row 2 shows the dynamics of the spatial averages of the solution components versus time (saplings component omitted).}\label{fig.waves_sup}
\end{figure} 
The simulations in Figure \ref{fig.waves_open} are for the same parameters as Figure 6 but in this case with an ``open'' boundary condition, as opposed to the reflective boundary condition used in the main text. These results show that the complex dynamics shown in Figure 6 are not dependent on the choice of boundary condition, but some small boundary effects are noticeable at the right-hand boundary in the third solution.
\begin{figure}[H]
	\centering
	\includegraphics[width=0.8\textwidth]{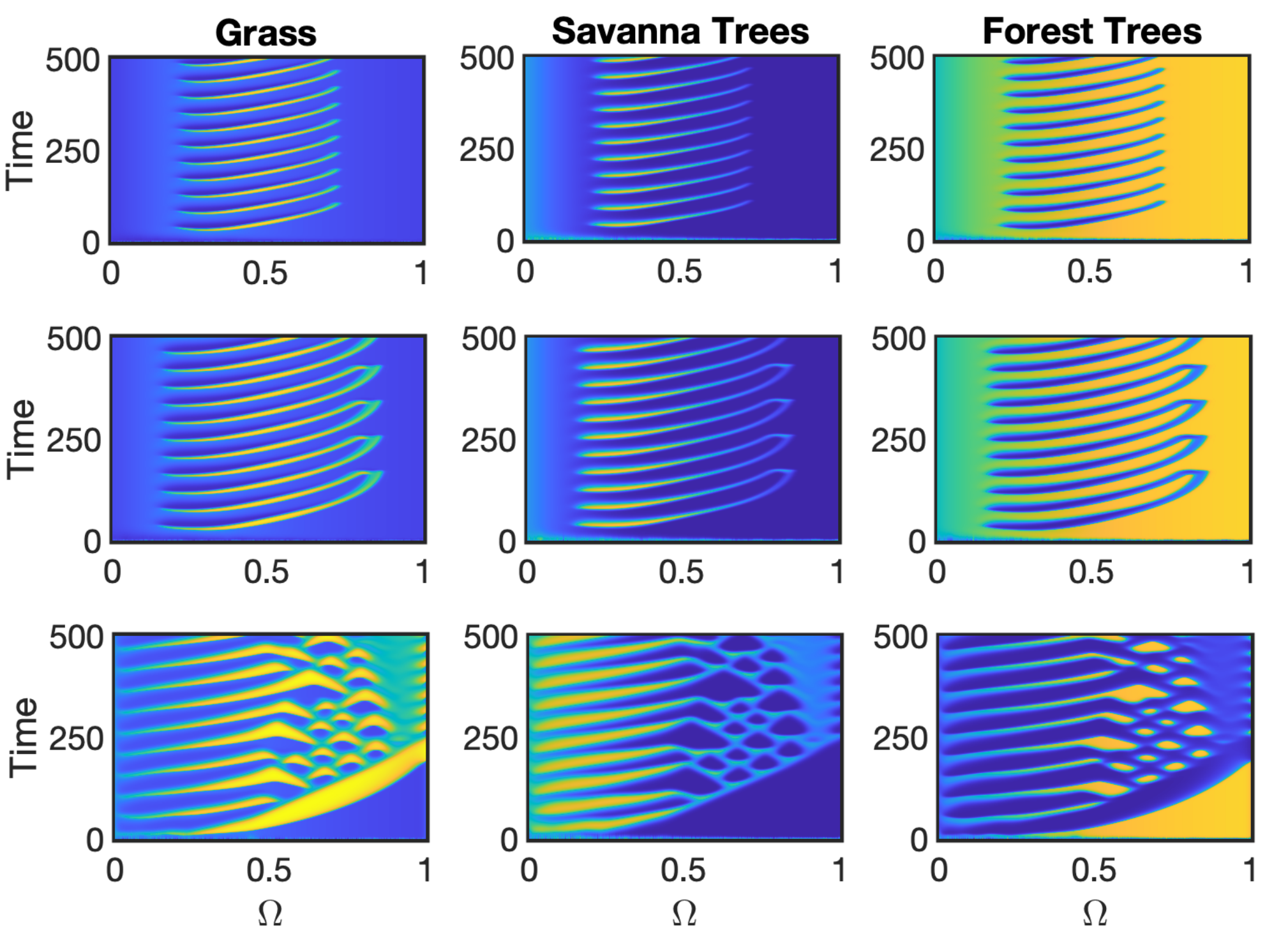}
	\caption{Solutions corresponding to those shown in Figure 6 for an ``open'' boundary condition.}\label{fig.waves_open}
\end{figure}
Figure \ref{fig.waves_summary_sup} shows a wider range of gradients intersecting the region of stable oscillations in $\alpha$-$\beta$ space. We vary the parameter $\beta_c$ from $0.2$ (top solution) to $2$ (bottom solution) and observe a transition from chaotic dynamics to increasing regular waves.
\begin{figure}[H]
	\centering
	\includegraphics[width=0.8\textwidth]{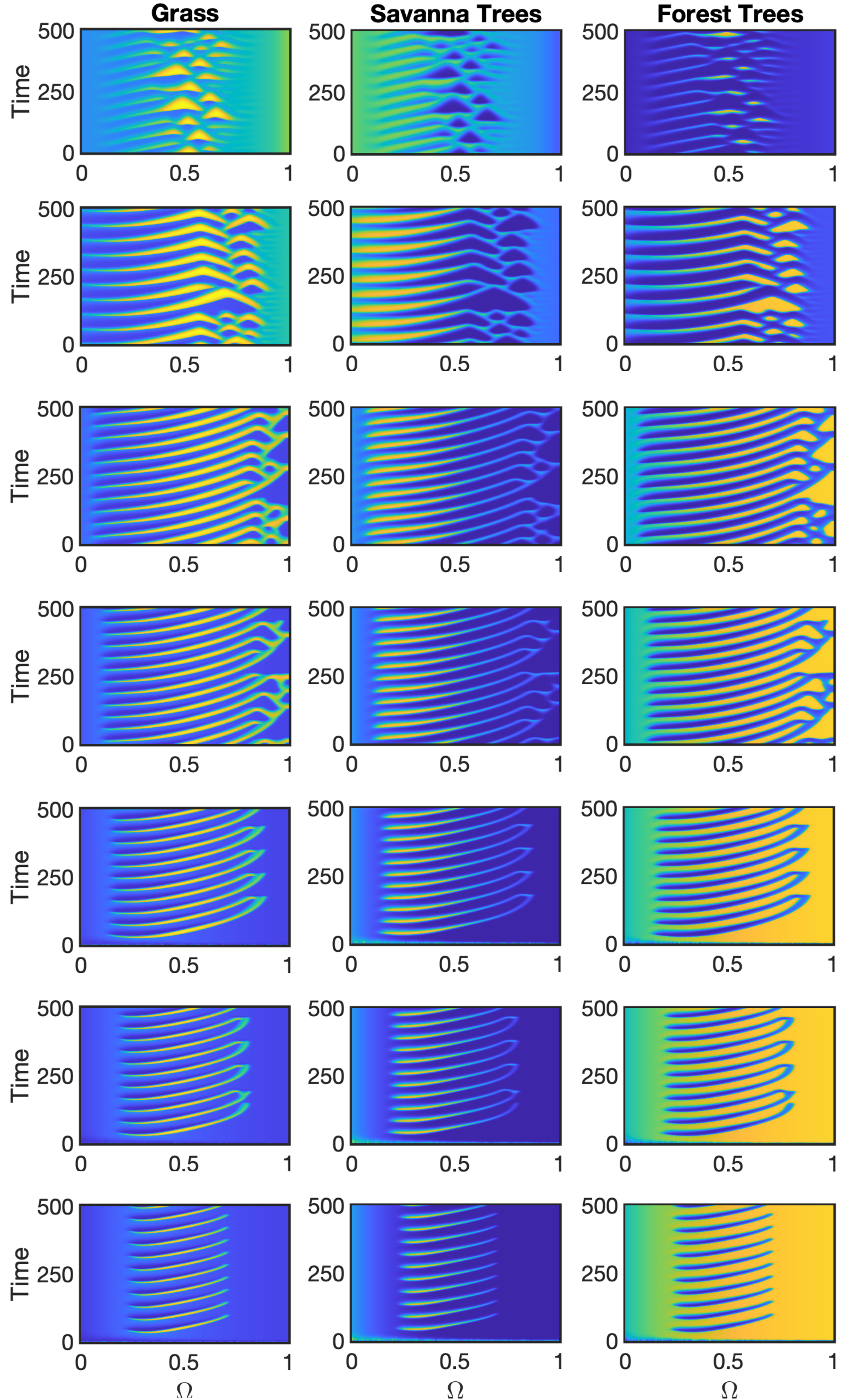}
	\caption{Solutions for a range of gradients between those shown in Figures 6A and 6C from the main text.}\label{fig.waves_summary_sup}
\end{figure}

\section{The Brain Patterning Model}\label{sec.brain_supp}
\subsection{Numerical parameters and schemes}
The PDEs were solved numerically using the open source finite element solver FreeFEM (version 4.2.1 - available at \href{https://freefem.org}{https://freefem.org}). We used a semi-implicit discretization scheme to evolve the weak formulation of the PDEs forward in time. We used a fixed step-size scheme, a uniformly spaced mesh and piecewise linear continuous finite elements ($P1$ elements in FreeFEM notation). For the 1D spatial simulations shown in the main text, we chose $\Omega = [0,40]$ and a uniformly spaced mesh with $400$ points, along with a stepsize of $h = 0.1$ (convergence was observed for $h<0.2$). The FreeFEM routines for solving the PDEs and the accompanying MATLAB code for processing output is available at \href{https://github.com/Touboul-Lab/cortex_patterning}{github.com/Touboul-Lab}.

\bibliographystyle{abbrv}
\bibliography{savanna_refs}

\begin{thebibliography}{10}

\bibitem{bastiaansen2020pulse}
R.~Bastiaansen, M.~Chirilus-Bruckner, and A.~Doelman.
\newblock Pulse solutions for an extended klausmeier model with spatially
  varying coefficients.
\newblock {\em SIAM Journal on Applied Dynamical Systems}, 19(1):1--57, 2020.

\bibitem{bastiaansen2022fragmented}
R.~Bastiaansen, H.~A. Dijkstra, and A.~S. von~der Heydt.
\newblock Fragmented tipping in a spatially heterogeneous world.
\newblock {\em Environmental Research Letters}, 17(4):045006, 2022.

\bibitem{belgacem1995effects}
F.~Belgacem and C.~Cosner.
\newblock The effects of dispersal along environmental gradients on the
  dynamics of populations in heterogeneous environment.
\newblock {\em Canadian Applied Mathematics Quarterly}, 3(4):379--397, 1995.

\bibitem{benson1998unravelling}
D.~L. Benson, P.~K. Maini, and J.~A. Sherratt.
\newblock Unravelling the {T}uring bifurcation using spatially varying
  diffusion coefficients.
\newblock {\em Journal of Mathematical Biology}, 37(5):381--417, 1998.

\bibitem{benson1993diffusion}
D.~L. Benson, J.~A. Sherratt, and P.~K. Maini.
\newblock Diffusion driven instability in an inhomogeneous domain.
\newblock {\em Bulletin of Mathematical Biology}, 55(2):365--384, 1993.

\bibitem{bucini2007continental}
G.~Bucini and N.~P. Hanan.
\newblock A continental-scale analysis of tree cover in {A}frican savannas.
\newblock {\em Global Ecology and Biogeography}, 16(5):593--605, 2007.

\bibitem{cantrell1996ecological}
R.~S. Cantrell, C.~Cosner, and V.~Hutson.
\newblock Ecological models, permanence and spatial heterogeneity.
\newblock {\em The Rocky Mountain Journal of Mathematics}, pages 1--35, 1996.

\bibitem{cantrell1993permanence}
R.~S. Cantrell, C.~Cosners, and V.~Hutson.
\newblock Permanence in ecological systems with spatial heterogeneity.
\newblock {\em Proceedings of the Royal Society of Edinburgh Section A:
  Mathematics}, 123(3):533--559, 1993.

\bibitem{champneys2021bistability}
A.~R. Champneys, F.~Al~Saadi, V.~F. Bre{\~n}a-Medina, V.~A. Grieneisen, A.~F.
  Mar{\'e}e, N.~Verschueren, and B.~Wuyts.
\newblock Bistability, wave pinning and localisation in natural
  reaction-diffusion systems.
\newblock {\em Physica D: Nonlinear Phenomena}, 416:132735, 2021.

\bibitem{crowley1989experimental}
M.~F. Crowley and I.~R. Epstein.
\newblock Experimental and theoretical studies of a coupled chemical
  oscillator: phase death, multistability and in-phase and out-of-phase
  entrainment.
\newblock {\em The Journal of Physical Chemistry}, 93(6):2496--2502, 1989.

\bibitem{dirr2006pinning}
N.~Dirr and N.~K. Yip.
\newblock Pinning and de-pinning phenomena in front propagation in
  heterogeneous media.
\newblock {\em Interfaces and Free Boundaries}, 8(1):79--109, 2006.

\bibitem{dornelas2021landscape}
V.~Dornelas, E.~H. Colombo, C.~L{\'o}pez, E.~Hern{\'a}ndez-Garc{\'\i}a, and
  C.~Anteneodo.
\newblock Landscape-induced spatial oscillations in population dynamics.
\newblock {\em Scientific Reports}, 11(1):1--11, 2021.

\bibitem{durrett2018heterogeneous}
R.~Durrett and R.~Ma.
\newblock A heterogeneous spatial model in which savanna and forest coexist in
  a stable equilibrium.
\newblock {\em arXiv preprint arXiv:1808.08159}, 2018.

\bibitem{durrett2015coexistence}
R.~Durrett and Y.~Zhang.
\newblock Coexistence of grass, saplings and trees in the {S}taver-{L}evin
  forest model.
\newblock {\em The Annals of Applied Probability}, 25(6):3434--3464, 2015.

\bibitem{ermentrout2003simulating}
B.~Ermentrout and A.~Mahajan.
\newblock Simulating, analyzing, and animating dynamical systems: a guide to
  {XPPAUT} for researchers and students.
\newblock {\em Applied Mechanics Reviews}, 56(4):B53--B53, 2003.

\bibitem{faye2014pulsatile}
G.~Faye and J.~Touboul.
\newblock Pulsatile localized dynamics in delayed neural field equations in
  arbitrary dimension.
\newblock {\em SIAM Journal on Applied Mathematics}, 74(5):1657--1690, 2014.

\bibitem{feng2021coup}
J.~Feng, W.-H. Hsu, D.~Patterson, C.-S. Tseng, H.-W. Hsing, Z.-H. Zhuang, Y.-T.
  Huang, A.~Faedo, J.~L. Rubenstein, J.~Touboul, and S.-J. Chou.
\newblock {COUP-TFI} specifies the medial entorhinal cortex identity and
  induces differential cell adhesion to determine the integrity of its boundary
  with neocortex.
\newblock {\em Science Advances}, 7(27):eabf6808, 2021.

\bibitem{flanagan:06}
J.~G. Flanagan.
\newblock Neural map specification by gradients.
\newblock {\em Current Opinion in Neurobiology}, 16(1):59--66, 2006.

\bibitem{garcia2000dispersal}
G.~Garc{\'\i}a-Ramos, F.~S{\'a}nchez-Gardu{\~n}o, and P.~K. Maini.
\newblock Dispersal can sharpen parapatric boundaries on a spatially varying
  environment.
\newblock {\em Ecology}, 81(3):749--760, 2000.

\bibitem{goel2020dispersal}
N.~Goel, V.~Guttal, S.~A. Levin, and A.~C. Staver.
\newblock Dispersal increases the resilience of tropical savanna and forest
  distributions.
\newblock {\em The American Naturalist}, 195(5):833--850, 2020.

\bibitem{goel2020dispersal_madagascar}
N.~Goel, E.~S. Van~Vleck, J.~C. Aleman, and A.~C. Staver.
\newblock Dispersal limitation and fire feedbacks maintain mesic savannas in
  {M}adagascar.
\newblock {\em Ecology}, 101(12):e03177, 2020.

\bibitem{he2018baptism}
T.~He and B.~B. Lamont.
\newblock Baptism by fire: the pivotal role of ancient conflagrations in
  evolution of the {E}arth's flora.
\newblock {\em National Science Review}, 5(2):237--254, 2018.

\bibitem{hillen2001global}
T.~Hillen and K.~Painter.
\newblock Global existence for a parabolic chemotaxis model with prevention of
  overcrowding.
\newblock {\em Advances in Applied Mathematics}, 26(4):280--301, 2001.

\bibitem{hillen2009user}
T.~Hillen and K.~J. Painter.
\newblock A user’s guide to {PDE} models for chemotaxis.
\newblock {\em Journal of Mathematical Biology}, 58(1):183--217, 2009.

\bibitem{hoyer2021impulsive}
A.~Hoyer-Leitzel and S.~Iams.
\newblock Impulsive fire disturbance in a savanna model: {T}ree--grass
  coexistence states, multiple stable system states, and resilience.
\newblock {\em Bulletin of Mathematical Biology}, 83(11):1--25, 2021.

\bibitem{hyman1991some}
B.~T. Hyman, G.~W. Van~Hoesen, and A.~R. Damasio.
\newblock Some cytoarchitectural abnormalities of the entorhinal cortex in
  schizophrenia.
\newblock {\em Archives of General Psychiatry}, 48(7):625--632, 1991.

\bibitem{keller1970initiation}
E.~F. Keller and L.~A. Segel.
\newblock Initiation of slime mold aggregation viewed as an instability.
\newblock {\em Journal of Theoretical Biology}, 26(3):399--415, 1970.

\bibitem{keller1971model}
E.~F. Keller and L.~A. Segel.
\newblock Model for chemotaxis.
\newblock {\em Journal of Theoretical Biology}, 30(2):225--234, 1971.

\bibitem{kiecker2005compartments}
C.~Kiecker and A.~Lumsden.
\newblock Compartments and their boundaries in vertebrate brain development.
\newblock {\em Nature Reviews Neuroscience}, 6(7):553--564, 2005.

\bibitem{kilpatrick2013wandering}
Z.~P. Kilpatrick and B.~Ermentrout.
\newblock Wandering bumps in stochastic neural fields.
\newblock {\em SIAM Journal on Applied Dynamical Systems}, 12(1):61--94, 2013.

\bibitem{kilpatrick2013optimizing}
Z.~P. Kilpatrick, B.~Ermentrout, and B.~Doiron.
\newblock Optimizing working memory with heterogeneity of recurrent cortical
  excitation.
\newblock {\em Journal of Neuroscience}, 33(48):18999--19011, 2013.

\bibitem{kilpatrick2008traveling}
Z.~P. Kilpatrick, S.~E. Folias, and P.~C. Bressloff.
\newblock Traveling pulses and wave propagation failure in inhomogeneous neural
  media.
\newblock {\em SIAM Journal on Applied Dynamical Systems}, 7(1):161--185, 2008.

\bibitem{klika2018domain}
V.~Klika, M.~Koz{\'a}k, and E.~A. Gaffney.
\newblock Domain size driven instability: {S}elf-organization in systems with
  advection.
\newblock {\em SIAM Journal on Applied Mathematics}, 78(5):2298--2322, 2018.

\bibitem{kozak2019pattern}
M.~Koz{\'a}k, E.~A. Gaffney, and V.~Klika.
\newblock Pattern formation in reaction-diffusion systems with piecewise
  kinetic modulation: {A}n example study of heterogeneous kinetics.
\newblock {\em Physical Review E}, 100(4):042220, 2019.

\bibitem{krause2021modern}
A.~L. Krause, E.~A. Gaffney, P.~K. Maini, and V.~Klika.
\newblock Modern perspectives on near-equilibrium analysis of {T}uring systems.
\newblock {\em Philosophical Transactions of the Royal Society A},
  379(2213):20200268, 2021.

\bibitem{krause2018heterogeneity}
A.~L. Krause, V.~Klika, T.~E. Woolley, and E.~A. Gaffney.
\newblock Heterogeneity induces spatiotemporal oscillations in
  reaction-diffusion systems.
\newblock {\em Physical Review E}, 97(5):052206, 2018.

\bibitem{krause2020one}
A.~L. Krause, V.~Klika, T.~E. Woolley, and E.~A. Gaffney.
\newblock From one pattern into another: {A}nalysis of {T}uring patterns in
  heterogeneous domains via {WKBJ}.
\newblock {\em Journal of the Royal Society Interface}, 17(162):20190621, 2020.

\bibitem{kulka1995influence}
A.~Kulka, M.~Bode, and H.-G. Purwins.
\newblock On the influence of inhomogeneities in a reaction-diffusion system.
\newblock {\em Physics Letters A}, 203(1):33--39, 1995.

\bibitem{li2019spatial}
Q.~Li, A.~C. Staver, W.~E, and S.~A. Levin.
\newblock Spatial feedbacks and the dynamics of savanna and forest.
\newblock {\em Theoretical Ecology}, 12(2):237--262, 2019.

\bibitem{mori2008wave}
Y.~Mori, A.~Jilkine, and L.~Edelstein-Keshet.
\newblock Wave-pinning and cell polarity from a bistable reaction-diffusion
  system.
\newblock {\em Biophysical Journal}, 94(9):3684--3697, 2008.

\bibitem{mori2011asymptotic}
Y.~Mori, A.~Jilkine, and L.~Edelstein-Keshet.
\newblock Asymptotic and bifurcation analysis of wave-pinning in a
  reaction-diffusion model for cell polarization.
\newblock {\em SIAM Journal on Applied Mathematics}, 71(4):1401--1427, 2011.

\bibitem{murray1982parameter}
J.~Murray.
\newblock Parameter space for {T}uring instability in reaction diffusion
  mechanisms: {A} comparison of models.
\newblock {\em Journal of Theoretical Biology}, 98(1):143--163, 1982.

\bibitem{nathan2012dispersal}
R.~Nathan, E.~Klein, J.~J. Robledo-Arnuncio, and E.~Revilla.
\newblock {\em Dispersal {K}ernels}, volume~15.
\newblock Oxford University Press, Oxford, UK, 2012.

\bibitem{nieto2018mesic}
P.~Nieto-Quintano, E.~T. Mitchard, R.~Odende, M.~A. Batsa~Mouwembe, T.~Rayden,
  and C.~M. Ryan.
\newblock The mesic savannas of the {B}ateke {P}lateau: {C}arbon stocks and
  floristic composition.
\newblock {\em Biotropica}, 50(6):868--880, 2018.

\bibitem{oleary:07}
D.~D. O'Leary, S.-J. Chou, and S.~Sahara.
\newblock Area patterning of the mammalian cortex.
\newblock {\em Neuron}, 56(2):252--269, 2007.

\bibitem{page2003pattern}
K.~Page, P.~K. Maini, and N.~A. Monk.
\newblock Pattern formation in spatially heterogeneous {T}uring
  reaction-diffusion models.
\newblock {\em Physica D: Nonlinear Phenomena}, 181(1-2):80--101, 2003.

\bibitem{patterson2019probabilistic}
D.~D. Patterson, S.~A. Levin, C.~Staver, and J.~D. Touboul.
\newblock Probabilistic foundations of spatial mean-field models in ecology and
  applications.
\newblock {\em SIAM Journal on Applied Dynamical Systems}, 19(4):2682--2719,
  2020.

\bibitem{pederick2018abnormal}
D.~T. Pederick, K.~L. Richards, S.~G. Piltz, R.~Kumar, S.~Mincheva-Tasheva,
  S.~A. Mandelstam, R.~C. Dale, I.~E. Scheffer, J.~Gecz, and S.~Petrou.
\newblock Abnormal cell sorting underlies the unique {X}-linked inheritance of
  {PCDH19} epilepsy.
\newblock {\em Neuron}, 97(1):59--66, 2018.

\bibitem{perthame2015competition}
B.~Perthame, C.~Qui{\~n}inao, and J.~Touboul.
\newblock Competition and boundary formation in heterogeneous media:
  {A}pplication to neuronal differentiation.
\newblock {\em Mathematical Models and Methods in Applied Sciences},
  25(13):2477--2502, 2015.

\bibitem{quesada2012basin}
C.~A. Quesada, O.~L. Phillips, M.~Schwarz, C.~I. Czimczik, T.~R. Baker,
  S.~Pati{\~n}o, N.~M. Fyllas, M.~G. Hodnett, R.~Herrera, and S.~Almeida.
\newblock Basin-wide variations in {A}mazon forest structure and function are
  mediated by both soils and climate.
\newblock {\em Biogeosciences}, 9(6):2203--2246, 2012.

\bibitem{schertzer2015implications}
E.~Schertzer, A.~Staver, and S.~Levin.
\newblock Implications of the spatial dynamics of fire spread for the
  bistability of savanna and forest.
\newblock {\em Journal of Mathematical Biology}, 70(1-2):329--341, 2015.

\bibitem{staal2020hysteresis}
A.~Staal, I.~Fetzer, L.~Wang-Erlandsson, J.~H. Bosmans, S.~C. Dekker, E.~H. van
  Nes, J.~Rockstr{\"o}m, and O.~A. Tuinenburg.
\newblock Hysteresis of tropical forests in the 21st century.
\newblock {\em Nature Communications}, 11(1):1--8, 2020.

\bibitem{staver2011tree}
A.~C. Staver, S.~Archibald, and S.~Levin.
\newblock Tree cover in sub-saharan africa: {R}ainfall and fire constrain
  forest and savanna as alternative stable states.
\newblock {\em Ecology}, 92(5):1063--1072, 2011.

\bibitem{staver_levin_2012}
A.~C. Staver and S.~Levin.
\newblock Integrating theoretical climate and fire effects on savanna and
  forest systems.
\newblock {\em The American Naturalist}, 180(2):211--224, 2012.

\bibitem{thompson2008plant}
S.~Thompson and G.~Katul.
\newblock Plant propagation fronts and wind dispersal: {A}n analytical model to
  upscale from seconds to decades using superstatistics.
\newblock {\em The American Naturalist}, 171(4):468--479, 2008.

\bibitem{touboul2018complex}
J.~D. Touboul, A.~C. Staver, and S.~A. Levin.
\newblock On the complex dynamics of savanna landscapes.
\newblock {\em Proceedings of the National Academy of Sciences},
  115(7):E1336--E1345, 2018.

\bibitem{van2015resilience}
I.~A. van~de Leemput, E.~H. van Nes, and M.~Scheffer.
\newblock Resilience of alternative states in spatially extended ecosystems.
\newblock {\em PLOS One}, 10(2), 2015.

\bibitem{van2021pattern}
R.~A. Van~Gorder.
\newblock Pattern formation from spatially heterogeneous reaction-diffusion
  systems.
\newblock {\em Philosophical Transactions of the Royal Society A},
  379(2213):20210001, 2021.

\bibitem{ward2002dynamics}
M.~J. Ward, D.~McInerney, P.~Houston, D.~Gavaghan, and P.~Maini.
\newblock The dynamics and pinning of a spike for a reaction-diffusion system.
\newblock {\em SIAM Journal on Applied Mathematics}, 62(4):1297--1328, 2002.

\bibitem{watrin2015causes}
F.~Watrin, J.-B. Manent, C.~Cardoso, and A.~Represa.
\newblock Causes and consequences of gray matter heterotopia.
\newblock {\em CNS Neuroscience \& Therapeutics}, 21(2):112--122, 2015.

\bibitem{wong2021spot}
T.~Wong and M.~J. Ward.
\newblock Spot patterns in the 2-{D} {S}chnakenberg model with localized
  heterogeneities.
\newblock {\em Studies in Applied Mathematics}, 146(4):779--833, 2021.

\bibitem{wuyts2017amazonian}
B.~Wuyts, A.~R. Champneys, and J.~I. House.
\newblock Amazonian forest-savanna bistability and human impact.
\newblock {\em Nature Communications}, 8(1):15519, 2017.

\bibitem{wuyts2019tropical}
B.~Wuyts, A.~R. Champneys, N.~Verschueren, and J.~I. House.
\newblock Tropical tree cover in a heterogeneous environment: A
  reaction-diffusion model.
\newblock {\em PLOS One}, 14(6), 2019.

\bibitem{xin2000front}
J.~Xin.
\newblock Front propagation in heterogeneous media.
\newblock {\em SIAM Review}, 42(2):161--230, 2000.

\end{thebibliography}
\end{document}